\documentclass[a4paper,10pt]{amsart}

\usepackage[english]{babel}
\usepackage{amsmath,amssymb,amsfonts,mathtools}
\usepackage{mathrsfs}
\usepackage{enumitem}
\usepackage{tikz}
\usetikzlibrary{arrows.meta}
\usepackage[expansion=false]{microtype}
\usepackage{hyperref}
\hypersetup{colorlinks=true,linkcolor=blue,citecolor=blue,urlcolor=blue}

\newtheorem{theorem}{Theorem}[section]
\newtheorem{proposition}[theorem]{Proposition}
\newtheorem{lemma}[theorem]{Lemma}
\newtheorem{corollary}[theorem]{Corollary}
\theoremstyle{definition}

\theoremstyle{remark}
\newtheorem{remark}[theorem]{Remark}

\newcommand{\R}{\mathbb R}
\newcommand{\N}{\mathbb N}
\newcommand{\Sph}{\mathbb S}
\newcommand{\eps}{\varepsilon}
\newcommand{\PV}{\operatorname{P.V.}}
\newcommand{\Ps}{\mathcal P_s}
\newcommand{\Psym}{\widehat{\mathcal P}_s}
\newcommand{\Bform}{\mathcal B_s}
\newcommand{\E}{\mathcal E}
\newcommand{\Wgt}{\mathcal W}

\newcommand{\dist}{\operatorname{dist}}
\newcommand{\Lout}{\mathrm{out}}
\newcommand{\supp}{\operatorname{supp}}
\newcommand{\Fou}{\mathcal F}

\DeclareMathOperator{\Rea}{Re}
\DeclareMathOperator{\Ima}{Im}

\title[Positive fractional bubble towers]{Positive fractional bubble towers for an equation with mixed Sobolev growth}

\author[S. Alarc\'on]{Salom\'on Alarc\'on}
\address{Departamento de Matem\'atica, Universidad T\'ecnica Federico Santa Mar\'ia, Av. Espa\~na 1680, Valpara\'iso, Chile}
\email{salomon.alarcon@usm.cl}

\author[L. M. Del Pezzo]{Leandro M. Del Pezzo}
\address{IESTA -- Departamento de M\'etodos Cuantitativos, Facultad de Ciencias Econ\'omicas y de Administraci\'on, Universidad de la Rep\'ublica, Av.~Gonzalo Ram\'irez 1926, 11200 Montevideo, Uruguay}
\email{leandro.delpezzo@fcea.edu.uy}

\author[A. Quaas]{Alexander Quaas}
\address{Departamento de Matem\'atica, Universidad T\'ecnica Federico Santa Mar\'ia, Av. Espa\~na 1680, Valpara\'iso, Chile}
\email{alexander.quaas@usm.cl}


\begin{document}

\begin{abstract}

	In this paper, we study the equation \((-\Delta)^s u=u^p+u^{p_s^*+\eps}\) in \(\R^N\), where \(0<s<1\), \(N\ge2\) and \(\frac{N}{N-2s}<p<p_s^*\coloneqq\frac{N+2s}{N-2s}\). For every integer \(n\ge1\) and every sufficiently small \(\eps>0\), we construct a bounded positive radial solution which, to leading order, is a tower of \(n\) critical bubbles centered at the origin. Moreover, we determine the concentration scales, whose limiting constants are given in terms of Beta functions, and we prove a uniform relative error estimate, the quantization of the critical mass and precise asymptotics at the origin and at infinity. As a consequence, the number of positive radial solutions becomes arbitrarily large as \(\eps\downarrow0\). This extends to the fractional setting the positive towers known for the local mixed-power equation. Our approach combines the Emden--Fowler transform, an explicit representation of the Green kernel of the cylindrical operator, weighted estimates for the projected linear problem and a finite-dimensional reduction. The proof covers the whole range \(0<s<1\) at once, including \(0<s\le\frac12\), where this Green kernel is unbounded and \(H^s(\R)\not\subset L^\infty(\R)\).

\end{abstract}

\maketitle

\section{Introduction}

	Concentration in critical elliptic problems is closely related to the dilation invariance of the limiting equation.
	The resulting loss of compactness is central to the fractional Yamabe and Brezis--Nirenberg problems.
	Profile decompositions describe this phenomenon through translations and dilations
	\cite{PalatucciPisante2014}, and global compactness on bounded domains was established in \cite{PalatucciPisante2015}.
	A bubble tower consists of critical profiles centered at one point, with asymptotically separated scales.
	We ask whether a small perturbation of the critical exponent produces positive entire solutions with any prescribed finite tower height.

	\medskip

	We consider the equation
	\[
		(-\Delta)^s u=u^p+u^q\qquad\text{in }\R^N,
	\]
	with
	\[
		u>0\quad\text{in }\R^N,
		\qquad u(x)\to0\quad\text{as }|x|\to\infty.
	\]
	We denote
	\[
		{p_s^*}\coloneqq\frac{N+2s}{N-2s},
		\qquad p_*^s\coloneqq\frac{N}{N-2s}=\frac{{p_s^*}+1}{2},
	\]
	and take
	\begin{equation}\label{eq:p-range}
		0<s<1,\qquad N\in\N,\quad N\ge2,
		\qquad p_*^s<p<{p_s^*}.
	\end{equation}

	We fix \(p\) and set
	\[
		q\coloneqq{p_s^*}+\eps,\qquad \eps\downarrow0.
	\]

	The fractional Laplacian is normalized according to
	\begin{equation}\label{eq:frac-lap}
		(-\Delta)^s u(x)
		=c_{N,s}\,\PV\int_{\R^N}
		\frac{u(x)-u(y)}{|x-y|^{N+2s}}\,dy,\qquad
		c_{N,s}=\frac{4^s s\,\Gamma(\frac N2+s)}
		{\pi^{N/2}\Gamma(1-s)}.
	\end{equation}
	This normalization gives the Fourier symbol \(|\zeta|^{2s}\); see \cite[Proposition~3.3]{DNPV} and \cite[(1.1.3)]{DipierroMedinaValdinoci2017}.

	\medskip

	The two powers lie on opposite sides of the critical exponent and allow a Pohozaev balance.
	For a single noncritical power, the corresponding identity excludes nontrivial solutions with finite energy and the required integrability.
	For the fast-decaying solutions constructed here, it reads
	\begin{equation}\label{eq:balance-intro}
		\frac{{p_s^*}-p}{p+1}\int_{\R^N}u^{p+1}\,dx
		=\frac{q-{p_s^*}}{q+1}\int_{\R^N}u^{q+1}\,dx.
	\end{equation}
	We derive this identity by translation in cylindrical variables, using the decay estimates proved below
	(Proposition~\ref{prop:balance}). For the bounded-domain identity, including its boundary term, see Ros-Oton and Serra
	\cite{RosOtonSerraPohozaev2014}.

	Both coefficients in \eqref{eq:balance-intro} are positive. Since \(p\) is fixed, the subcritical integral is
	\(O(\eps)\) relative to the supercritical one. Concentration produces this separation.

	\medskip

	The lower bound on \(p\) ensures the integrability of the subcritical term. For a bounded solution with
	\[
		C^{-1}|x|^{2s-N}\le u(x)\le C|x|^{2s-N}
		\qquad(|x|\ge R)
	\]
	for some \(C\ge1\) and \(R>0\), we have
	\[
		u^p\in L^1(\R^N\setminus B_1)
		\quad\text{precisely when}\quad p(N-2s)>N.
	\]
	The condition \(p>p_*^s\) is also necessary for positive entire solutions: since \((-\Delta)^su\ge u^p\), the Liouville theorem for supersolutions \cite[Theorem~1.3]{FelmerQuaas2011}, specialized to \((-\Delta)^s\), excludes positive solutions when \(1<p\le p_*^s\).
	In the cylindrical variables introduced in Section~\ref{sec:EF}, the inequality \(p>p_*^s\) becomes
	\begin{equation}\label{eq:2pP}
		2p-{p_s^*}>1.
	\end{equation}
	Thus the subcritical term decays faster than the Green kernel as \(t\to-\infty\), that is, as \(|x|\to\infty\). At equality, the rates coincide and the spatial integral diverges.

	\medskip

	Set
	\[
		a\coloneqq\frac{N-2s}{2},\qquad b\coloneqq{p_s^*}-p>0.
	\]
	The critical equation is
	\[
		(-\Delta)^s W=W^{p_s^*}\qquad\text{in }\R^N.
	\]
	Its standard bubbles are
	\[
		W_{\mu,\zeta}(x)
		=\alpha_{N,s}\left(\frac{\mu}{\mu^2+|x-\zeta|^2}\right)^a,
		\qquad \mu>0,\quad\zeta\in\R^N,
	\]
	where
	\begin{equation}\label{eq:alpha-value}
		\alpha_{N,s}^{{p_s^*}-1}
		=4^s\frac{\Gamma(\frac{N+2s}{2})}{\Gamma(\frac{N-2s}{2})}.
	\end{equation}
	With this normalization, the standard profile satisfies
	\[
		(-\Delta)^s(1+|x|^2)^{-a}
		=\alpha_{N,s}^{{p_s^*}-1}
		(1+|x|^2)^{-a-2s}, \qquad a p_s^*=a+2s;
	\]
	see \cite[(1.5)--(1.6)]{DavilaDelPinoSire2013}. The value of \(\alpha_{N,s}\) follows from the Riesz potential identity \cite[(3.10)--(3.13)]{MR0717827}. Lieb's sharp Hardy--Littlewood--Sobolev theorem identifies the corresponding extremals \cite[Theorem~3.1]{MR0717827}. The classification results of Li and Chen--Li--Ou \cite{MR2055032,ChenLiOu2006}
	characterize the positive critical profiles under their regularity and integrability assumptions. Together with the profile decomposition, these results suggest an approximate solution formed by critical bubbles at separated scales. The reduction determines the scales for which this superposition can be corrected to an exact solution.

\subsection*{Related results}
	For the Laplacian, Lin and Ni \cite[Proposition~2.8]{MR0920985} obtained explicit slowly decaying solutions when \(q=2p-1\). Bam\'on, Flores and 
	del Pino \cite{MR1791878} later studied radial ground states through a geometric analysis of the Emden--Fowler system. For fixed 
	\(\frac N{N-2}<p<\frac{N+2}{N-2}\), they proved that the number of radial ground states with fast decay becomes arbitrarily large as 
	\(q\downarrow\frac{N+2}{N-2}\) \cite[Theorem~1.1(b)]{MR1791878}. Corollary~\ref{cor:multiplicity} is the fractional counterpart of this result. 
	Bam\'on, Flores and del Pino also proved that, for fixed \(q\), there are no radial ground states when \(p\) is close to \(\frac N{N-2}\) or below it 
	\cite[Theorem~1.2]{MR1791878}. Flores \cite{Flores2004} studied the resonance between the two powers at \(q=2p-1\) and its 
	consequences for multiplicity. Since \(p_*^s=\frac{p_s^*+1}2\), the integrability threshold \(p=p_*^s\) considered here is exactly where the line \(q=2p-1\), resonant in the local case, meets \(q=p_s^*\). D\'avila and Guerra \cite{MR3540603} 
	constructed slowly decaying solutions when \(p\) is close to \(\frac N{N-2}\) from above, \(q\) is close to \(\frac{N+2}{N-2}\) and a 
	relation between them holds.

	\medskip

	The closest local analogue is due to Campos \cite{MR2381679}, who constructed, for fixed \(\frac N{N-2}<p<\frac{N+2}{N-2}\) and \(q=\frac{N+2}{N-2}+\eps\), positive radial towers of arbitrary prescribed height and determined their leading scales. The logarithmic parametrization of the centers and the two balances used below are the fractional counterparts of that construction. 

	\medskip

	For equations with a mass term, D\'avila, del Pino and Guerra \cite{DavilaDelPinoGuerra2013} proved nonuniqueness of positive ground states. 
	Chen, D\'avila and Guerra \cite{ChenDavilaGuerra2016} constructed positive towers for
	\[
		-\Delta u+u=u^{\frac{N+2}{N-2}+\eps}+\lambda u^\nu
		\qquad\text{in }\R^N,
	\]
	with \(\lambda>0\) and a fixed subcritical exponent \(\nu\) in the range of their result. The mass term changes the decay at infinity. Finite-dimensional reduction has also been used for slightly supercritical Brezis--Nirenberg problems, in the ball \cite{dPDM2003} and in general domains \cite{GeJingPacard2005,MussoPistoia2010}, and for nodal radial solutions of slightly subcritical problems in the ball \cite{ContrerasDelPino2006}.

	\medskip

	The fractional Brezis--Nirenberg problem on bounded domains was studied in \cite{ServadeiValdinoci2015}, and the concentration of least-energy solutions of nearly critical problems in \cite{ChoiKimLee2014}; whole-space critical problems were treated in \cite{DipierroMedinaValdinoci2017}. Positive solutions and concentrating standing waves for fractional Schr\"odinger equations were obtained in \cite{FelmerQuaasTan2012,DavilaDelPinoWei2014}, and uniqueness and nondegeneracy of ground states in \cite{FrankLenzmannSilvestre2016}.

	\medskip

	Closer to the present construction, D\'avila, L\'opez R\'ios and Sire \cite{DavilaLopezRiosSire2017} developed the finite-dimensional reduction for the 
	slightly subcritical and supercritical problems \((-\Delta)^su=u^{{p_s^*}\pm\eps}\) in bounded domains, for the spectral and the restricted fractional 
	Laplacian. Their solutions concentrate at distinct points, which are determined by the Green and Robin functions of the domain.
	In \(\R^N\), Wang and Wei \cite{WangWei2021} used a finite-dimensional reduction and local Pohozaev identities to construct, under their hypotheses on \(N,s\) and the nonnegative bounded potential \(K\), solutions of \((-\Delta)^su=K(x)u^{p_s^*\pm\eps}\) with a number of bubbles that tends to infinity as \(\eps\downarrow0\). These bubbles are centered at the vertices of a regular polygon, whereas here all bubbles are centered at the origin and the multiplicity comes from towers of different heights.

	\medskip

	Positive fractional towers have also been constructed in geometric problems.
	Chen and Kim \cite{ChenKim2018} obtained towers for asymptotically critical perturbations of the fractional Yamabe equation on conformal infinities; see also \cite{AoDelaTorreGonzalezWei2020}.
	DelaTorre, del Pino, Gonz\'alez and Wei \cite{MR3694655} studied Delaunay-type singular
	solutions. Here we construct regular entire solutions of an autonomous mixed-power equation.

	\medskip

	For nodal fractional towers, Cora and Iacopetti \cite{CoraIacopetti2018} studied least-energy radial sign-changing solutions of the fractional Brezis--Nirenberg problem in a ball, for \(N>6s\). For \(\frac12<s<1\), if these solutions \(u_\lambda\) change sign exactly once for all sufficiently small \(\lambda>0\) and are normalized to be nonnegative near the origin, then, as \(\lambda\downarrow0\), their positive and negative parts concentrate at the origin with different speeds, and the rescaled positive part converges to the standard bubble \cite[Theorem~1.4]{CoraIacopetti2018}. For the slightly subcritical problem \((-\Delta)^su=|u|^{{p_s^*}-1-\eps}u\) in a ball, with \(0<s<1\) and \(N>2s\), they showed that, as \(\eps\downarrow0\), the positive and negative parts of least-energy radial sign-changing solutions \(u_\eps\), normalized by \(u_\eps(0)>0\), blow up with different speeds, the positive part concentrates at the center and, after rescaling, converges to the standard bubble; for \(\frac12<s<1\), the negative part also concentrates at the center \cite[Theorem~1.1]{CoraIacopetti2019}.

	\medskip

	Chen, Long and Yang \cite{ChenLongYang2022} constructed sign-changing towers for a critical problem in domains with a small hole. More recently, Yu and Zhang \cite{YuZhang2026} constructed sign-changing towers for an almost critical non-power problem involving the spectral fractional Laplacian in bounded domains. Cai, Wang, Wei and Yang \cite{CaiWangWeiYang2025} obtained a dynamical construction for the critical fractional heat equation: positive ancient towers and sign-changing towers as forward time tends to infinity. The towers constructed here are stationary and positive.

\subsection*{The scales}
	Two balances determine the concentration scales. To explain them, we write
	\[
		\xi_i\coloneqq-a\log\mu_i,
		\qquad \xi_1<\cdots<\xi_n.
	\]
	The subcritical term is of order \(e^{-b\xi_1}\), whereas the interaction between consecutive bubbles is of order \(e^{-(\xi_i-\xi_{i-1})}\). The perturbation of the upper exponent contributes at order \(\eps\) to the reduced equations. This suggests the balances
	\[
		C^{-1}\eps\le e^{-b\xi_1}\le C\eps,\qquad
		C^{-1}\eps\le e^{-(\xi_i-\xi_{i-1})}\le C\eps
		\qquad(2\le i\le n).
	\]
	These balances fix the largest scale and the successive scale ratios; the reduced energy gives their leading constants.

\subsection*{The main result}
	The positive constants \(a_2,a_3,a_5\) are given by single-bubble integrals in \eqref{eq:a3}--\eqref{eq:a2} and evaluated in Remark~\ref{rem:explicit}. The H\"older spaces \(C^{2s+\sigma}\) follow the convention of Section~\ref{sec:EF}.
	The threshold \(\eps_n\) and the constants in the estimates may depend on \(N,s,p,n,\gamma\) and \(\alpha\).

	\begin{theorem}\label{thm:main}
		Assume \eqref{eq:p-range}, fix \(n\in\N\) and choose
		\begin{equation}\label{eq:gamma-alpha-thm}
			0<\gamma<\min\{1,{p_s^*}-1\},
			\qquad
			\frac12<\alpha<\min\Big\{1,\frac{{p_s^*}-\gamma}{2}\Big\}.
		\end{equation}
		There exist \(\eps_n>0\) and \(C>0\) such that, for every \(0<\eps<\eps_n\), the problem
		\begin{equation}\label{eq:main-problem}
			\begin{cases}
				(-\Delta)^s u=u^p+u^{{p_s^*}+\eps}&\text{in }\R^N,\\
				u>0&\text{in }\R^N,\\
				u(x)\to0&\text{as }|x|\to\infty
			\end{cases}
		\end{equation}
		has a radial solution \(u_{\eps,n}\in L^\infty(\R^N)\), which belongs to
		\[
			C^{2s+\sigma}(\R^N)
			\qquad\text{for every }0<\sigma<\min\{1,2-2s\}\text{ with }2s+\sigma\ne1,
		\]
		and has the following properties.
		\begin{enumerate}[label=\textup{(\roman*)}]
			\item \emph{Scales.} There are parameters
			\[
				\Lambda_\eps=(\Lambda_{1,\eps},\ldots,\Lambda_{n,\eps})
				\to\Lambda^*\in(0,\infty)^n
			\]
			and scales satisfying
			\begin{align}
				&\mu_{1,\eps}
				=\eps^{1/(ab)}\Lambda_{1,\eps}^{1/a},\label{eq:mu1}\\
				&\frac{\mu_{i,\eps}}{\mu_{i-1,\eps}}
				=(\eps\Lambda_{i,\eps})^{1/a},
				\qquad 2\le i\le n.\label{eq:muratio}
			\end{align}
			In particular,
			\[
				\mu_{1,\eps}\gg\mu_{2,\eps}\gg\cdots\gg\mu_{n,\eps}\to0.
			\]
			\item \emph{Tower decomposition.} The solution has the form
			\[
				u_{\eps,n}(x)
				=\sum_{i=1}^nW_{\mu_{i,\eps},0}(x)+\omega_{\eps,n}(x).
			\]
			The remainder \(\omega_{\eps,n}\) is continuous on \(\R^N\). For \(x\ne0\),
			\begin{equation}\label{eq:remainder-x}
				\begin{split}
					|\omega_{\eps,n}(x)|
					\le C\eps^\alpha\sum_{i=1}^n\mu_{i,\eps}^{-a}
					\min\biggl\{&\left(\frac{|x|}{\mu_{i,\eps}}\right)^{-a(1-\gamma)},\left(\frac{|x|}{\mu_{i,\eps}}\right)^{-a(1+\gamma)}\biggr\}.
				\end{split}
			\end{equation}
			\item \emph{Fast decay.} The largest scale satisfies
			\[
				\mu_{1,\eps}^a=e^{-\xi_{1,\eps}}
				=\eps^{1/b}\Lambda_{1,\eps}.
			\]
			For \(|x|\ge\mu_{1,\eps}\),
			\begin{equation}\label{eq:fast-decay}
				C^{-1}\mu_{1,\eps}^{a}|x|^{2s-N}
				\le u_{\eps,n}(x)
				\le C\mu_{1,\eps}^{a}|x|^{2s-N}.
			\end{equation}
			\item \emph{Limiting parameters.} The limiting values are
			\begin{equation}\label{eq:Lambda-star-intro}
				\Lambda_1^*=\Big(\frac{na_3}{ba_5}\Big)^{1/b},
				\qquad
				\Lambda_i^*=\frac{a_3(n-i+1)}{a_2},\quad 2\le i\le n.
			\end{equation}
			Consequently, with \(D_i=(\prod_{j=1}^i\Lambda_j^*)^{1/a}\),
			\begin{equation}\label{eq:mu-asymp-intro}
				\mu_{i,\eps}
				=D_i\eps^{\frac1a(\frac1b+i-1)}(1+o(1)),
				\qquad 1\le i\le n.
			\end{equation}
		\end{enumerate}
	\end{theorem}

	\begin{remark}
		\begin{enumerate}[label=\textup{(\alph*)}]
			\item The parameters \(\Lambda_i^*\) for \(i\ge2\) are independent of \(p\) and unchanged if the fractional Laplacian is multiplied by a positive constant. The first parameter depends on the bubble amplitude \(\alpha_{N,s}\) and hence on this normalization.
			\item In cylindrical variables (Sections~\ref{sec:energy} and~\ref{sec:proof-main}), the solution has energy
			\[
				\E_\eps(v_{\eps,n})=na_1+O(\eps|\log\eps|),
			\]
			where \(a_1>0\) is the energy of one critical bubble.
			\item The condition \(\alpha>\frac12\) ensures that the energy contribution of the correction is \(o(\eps)\). The upper bound \(\alpha<\frac{{p_s^*}-\gamma}{2}\) comes from the weighted estimates for the interactions.
			\item The solution is obtained from a maximizer \(\Lambda_\eps\) of the reduced energy and a fixed point in a ball of radius \(C_0\eps^\alpha\). We do not claim that it is independent of the choice of \((\gamma,\alpha)\) or of the maximizer; all the asymptotic statements hold for every such choice.
			\item The auxiliary parameters \(\eta\) and \(\delta\) are fixed in terms of \(N,s,p,n\) in Sections~\ref{sec:EF} and~\ref{sec:proof-main}. The estimates are established for fixed \(n\), are not uniform near the endpoints of \eqref{eq:p-range} and give no quantitative regime with \(n\to\infty\).
		\end{enumerate}
	\end{remark}

	The exterior estimates yield a global relative approximation and precise asymptotics at the origin and at infinity.

	\begin{corollary}\label{cor:asymptotics}
		For the solutions of Theorem~\ref{thm:main}, set
		\[
			\varrho\coloneqq\alpha-\frac{1-\gamma}{2}>0.
		\]
		Then the following statements hold as \(\eps\downarrow0\).
		\begin{enumerate}[label=\textup{(\roman*)}]
			\item The relative error satisfies
			\[
				\sup_{x\in\R^N}
				\frac{\bigg|u_{\eps,n}(x)-\displaystyle\sum_{i=1}^nW_{\mu_{i,\eps},0}(x)\bigg|}
				{\displaystyle\sum_{i=1}^nW_{\mu_{i,\eps},0}(x)}
				\le C\eps^\varrho.
			\]
			The quotient at \(x=0\) is defined by continuity. For \(n=1\), the right-hand side can be replaced by \(C\eps^\alpha\).
			\item The value at the origin satisfies
			\[
					u_{\eps,n}(0)
					=\alpha_{N,s}\sum_{i=1}^n\mu_{i,\eps}^{-a}
					+O(\eps^\alpha\mu_{n,\eps}^{-a})=\alpha_{N,s}\mu_{n,\eps}^{-a}(1+o(1)).
			\]
			\item For each \(1\le i\le n\),
			\[
				\mu_{i,\eps}^a u_{\eps,n}(\mu_{i,\eps}y)
				\to W_{1,0}(y)
			\]
			uniformly on compact subsets of \(\R^N\setminus\{0\}\). For \(i=n\), convergence is uniform on compact subsets of \(\R^N\).
			\item The critical mass is quantized:
			\[
				\int_{\R^N}u_{\eps,n}^{{p_s^*}+1}\,dx
				\to n\int_{\R^N}W_{1,0}^{{p_s^*}+1}\,dx.
			\]
			\item The following limit exists and is positive:
			\[
					A_{\eps,n}
					\coloneqq\lim_{|x|\to\infty}|x|^{N-2s}u_{\eps,n}(x)
					=\alpha_{N,s}\sum_{i=1}^n\mu_{i,\eps}^{a}
					+O(\eps^\alpha\mu_{1,\eps}^{a})
					=\alpha_{N,s}\mu_{1,\eps}^{a}(1+o(1)).
			\]
		\end{enumerate}
	\end{corollary}

	\begin{corollary}\label{cor:multiplicity}
		For every integer \(m\ge1\), there exists \(\eps_m'>0\) such that \eqref{eq:main-problem} has at least \(m\) distinct positive radial solutions whenever \(0<\eps<\eps_m'\). Since \(\eps_m'\) depends on \(m\), this does not give infinitely many solutions for a fixed \(\eps>0\).
	\end{corollary}

	\begin{remark}
		The following questions are not treated here.
		\begin{enumerate}[label=\textup{(\alph*)}]
			\item Other parameter regimes and slowly decaying solutions. The construction concerns fixed \(p_*^s<p<p_s^*\) and \(q\downarrow p_s^*\).
			\item Local uniqueness, the radial Morse index and nonradial nondegeneracy modulo translations.
			\item The case \(N=1\), admissible only for \(0<s<\frac12\). The sphere \(\Sph^{N-1}\) is replaced by \(\Sph^0=\{\pm1\}\), radial functions are even and the kernel analysis requires a separate treatment.
		\end{enumerate}
	\end{remark}

\subsection*{Main contributions}
	The main features of the construction are the following.
	\begin{enumerate}[label=\textup{(\alph*)}]
		\item For each fixed \(0<s<1\), each \(p_*^s<p<p_s^*\) and each prescribed height \(n\ge1\), we construct positive radial towers with explicit leading scales and a uniform relative error in \(\R^N\).
		This extends Campos' local mixed-power construction \cite{MR2381679} to the fractional setting.
		\item We obtain an explicit formula for the Green kernel \(G_s\) of the cylindrical operator \(\Ps\). Its positive Riesz potential representation gives the exact coefficient \(\beta_s^{-1}\) of its exponential tail, which determines the leading term of the solutions at infinity, and its 
		precise behavior at the origin: near \(0\), \(G_s\) is bounded for \(\frac12<s<1\) and comparable to \(\log\frac1{|z|}\) for \(s=\frac12\) 
		and to \(|z|^{2s-1}\) for \(0<s<\frac12\) (Lemma~\ref{lem:indicial}). For the symbol of \(\Ps\), see \cite{AoChanDelaTorreEtAl2019}. The representation of the singular kernel \(k_s\) and its behavior at \(0\) and at infinity are given in \cite[(2.5)--(2.6) and Lemma~2.6]{MR3694655}, in the variable \(\log r\).
		\item The same reduction applies for each fixed \(0<s<1\). In the range \(0<s\le\frac12\), the embedding \(H^s(\R)\subset L^\infty(\R)\) 
		fails. Hence the solutions of the projected problem obtained in the energy space \(H^s_k(\R)\) need not be bounded a priori, and the Green kernel 
		is unbounded. Boundedness, decay and H\"older continuity of the solutions follow from the Green 
		representation, Young's inequality and the comparison of \(G_s\) with \(e^{-|z|}\) on slowly varying majorants. Higher regularity follows from the Schauder estimates of 
		\cite{MR3694655,DongKim2013} on the line and of \cite{SilvestreCPAM2007} in \(\R^N\). The same regularity bootstrap supplies the global bound on \(v'\) used in the translation argument for the Pohozaev balance.
	\end{enumerate}

\subsection*{Main difficulties and outline of the proof}
	The change of variables
	\[
		t=-a\log r,\qquad u(r)=r^{-a}v(t)
	\]
	turns dilations into translations. In these variables, the radial equation becomes
	\[
		\Ps v=e^{-bt}v^p+e^{\eps t}v^{p_s^*+\eps}.
	\]
	Here \(\Ps\) is the rescaled radial conformal fractional Laplacian on the cylinder.

	\medskip

	The nonlocal problem requires particular care at four points in the reduction.
	\begin{enumerate}[label=\textup{(\roman*)}]
		\item \emph{The projected linear problem.} The local phase-plane argument is unavailable. We use compactness for the nonlocal operator and single-bubble nondegeneracy \cite[Theorem~1.1]{DavilaDelPinoSire2013}. To apply the latter, we first prove exponential decay for each bounded cylindrical kernel element and removability of the corresponding singularity at the origin.
		\item \emph{Resonance.} The bubble, its translation mode and the Green kernel all decay like \(e^{-|t|}\). Convolution at this rate produces the factor \((1+|t|)e^{-|t|}\). The exterior norm includes this factor while the projection coefficients are present. Once they vanish, a bootstrap argument gives the sharper decay.
		\item \emph{The weight in the subcritical term.} The coefficient \(e^{-bt}\) grows as \(t\to-\infty\) and cannot be controlled by the weaker exponential norm alone. We therefore estimate the right-hand side in a second weighted norm with exterior decay exponent \(\eta>1\). The inequality \(2p-p_s^*>1\) ensures that the nonlinear term belongs to this space.
		\item \emph{Positivity and the reduced energy.} We use the positive-part nonlinearity in the equation and the energy, and recover positivity from the relative error estimate.
		The reduced energy is \(C^1\), but only its uniform asymptotic expansion is needed: a strict limiting maximum determines the scales.
		No expansion of its derivatives is required.
	\end{enumerate}

	\medskip

	The paper is organized as follows. In Section~\ref{sec:EF}, we develop the cylindrical formulation, prove the radial nondegeneracy of the bubble and introduce the weighted spaces. 
	Sections~\ref{sec:linear} and~\ref{sec:nonlinear} are devoted to the projected linear and nonlinear problems, respectively. In Section~\ref{sec:energy}, we compute the reduced energy. In 
	Section~\ref{sec:proof-main}, we determine the scales and prove the asymptotic formulas and the Pohozaev balance. Finally, in Appendix~\ref{app:functional}, 
	we collect the functional and regularity estimates.

\section{Emden--Fowler formulation and the critical profile}\label{sec:EF}

	We write the radial problem on \(\R\times\Sph^{N-1}\). For \(0<\varkappa<2\), set \(C^\varkappa=C^{0,\varkappa}\) if \(\varkappa<1\), \(C^\varkappa=C^1\) if \(\varkappa=1\), and \(C^\varkappa=C^{1,\varkappa-1}\) if \(\varkappa>1\), with the same convention for local spaces. Set
	\[
		\varpi_s(z)\coloneqq
		\begin{cases}
			1&\text{if }\frac12<s<1,\\
			1+\log\frac1{|z|}&\text{if }s=\frac12,\\
			|z|^{2s-1}&\text{if }0<s<\frac12,
		\end{cases}
		\qquad 0<|z|\le1,
	\]
	and \(\varsigma_s\coloneqq\infty\) if \(s\ge\frac12\), \(\varsigma_s\coloneqq1/(1-2s)\) if \(s<\frac12\). Then \(\varpi_s\in L^{\varsigma}(-1,1)\) for every \(1\le\varsigma<\varsigma_s\).

	\medskip

	For \(0<s\le\frac12\), the Green kernel in Lemma~\ref{lem:indicial} is singular at the origin and comparable to \(\varpi_s\), and \(H^s(\R)\not\subset L^\infty(\R)\). The same proof applies throughout \(0<s<1\), using:
	\begin{enumerate}[label=\textup{(\alph*)}]
		\item integrability of \(G_s\) and its tail \(\beta_s^{-1}e^{-|z|}\), which give comparison with \(e^{-|z|}\) on slowly varying majorants (Lemma~\ref{lem:comparison});
		\item bounded representatives in the energy space, obtained from the Green representation and Young's inequality (Lemma~\ref{lem:bounded-rep});
		\item H\"older continuity from the Green representation and the bound on \(G_s'\) in Lemma~\ref{lem:indicial}(iii), followed by a finite Schauder bootstrap (Proposition~\ref{prop:regularity}).
	\end{enumerate}

\subsection{The transformed operator}

	Let \(u=u(r)\) be radial and write
	\[
		u(r)=r^{-a}v(t)=e^t v(t),
		\qquad
		t=-a\log r
		.
	\]
	Thus \(r=e^{-t/a}\) and
	\[
		r\to\infty\iff t\to-\infty,
		\qquad
		r\to0\iff t\to+\infty.
	\]
	We use the identities
	\[
		a{p_s^*}=a+2s=N-a=\frac{N+2s}{2},
		\qquad
		a({p_s^*}+1)=N,
		\qquad
		\frac{N+2s}{a}=2{p_s^*} .
	\]
	For \(\rho>0\), \(\rho\neq1\), set
	\[
		K(\rho)\coloneqq\int_{\Sph^{N-1}}\frac{d\omega}{|e_1-\rho\,\omega|^{N+2s}},
		\qquad e_1=(1,0,\dots,0).
	\]
	The kernel \(K\) is analytic away from \(\rho=1\). Since
	\[
		|e_1-\rho\omega|^2=1-2\rho\omega_1+\rho^2=|\rho e_1-\omega|^2,
	\]
	it satisfies
	\begin{equation}\label{eq:K-inversion}
		K(1/\rho)=\rho^{N+2s}K(\rho).
	\end{equation}

	\begin{proposition}[Radial Emden--Fowler representation]\label{prop:EF}
		Define, for \(\tau\neq0\),
		\[
			k_s(\tau)\coloneqq\frac{c_{N,s}}{a}\,e^{{p_s^*}\tau}K\bigl(e^{\tau/a}\bigr),
			\qquad
			m_s\coloneqq\int_{\R}(\cosh\tau-1)\,k_s(\tau)\,d\tau .
		\]
		\begin{enumerate}[label=\textup{(\roman*)}]
			\item \(k_s\) is even, positive and real-analytic on \(\R\setminus\{0\}\), and there are constants \(c_1,c_2>0\) and \(C>0\) such that
			\begin{equation}\label{eq:kernel-expansion}
				k_s(\tau)=c_1|\tau|^{-1-2s}+\tilde k(\tau),
				\qquad |\tilde k(\tau)|\le C\left(|\tau|^{1-2s}+\log\frac2{|\tau|}\right)\quad(0<|\tau|\le1),
			\end{equation}
			\begin{equation}\label{eq:kernel-infinity}
				k_s(\tau)=c_2e^{-{p_s^*}|\tau|}\bigl(1+O(e^{-2|\tau|/a})\bigr)\qquad(|\tau|\to\infty),
			\end{equation}
			with \(c_1=a^{2s}c_{1,s}\) and \(c_2=c_{N,s}|\Sph^{N-1}|/a\), where \(c_{1,s}\) is the normalization constant in \eqref{eq:frac-lap} in dimension one.  In particular, \(\tilde k\in L^1(\R)\) (with \(\tilde k\coloneqq k_s-c_1|\tau|^{-1-2s}\) on all of \(\R\setminus\{0\}\)), \(0<m_s<\infty\), and
			\begin{equation}\label{eq:kernel-bounds}
				\begin{aligned}
					C^{-1}|\tau|^{-1-2s}&\le k_s(\tau)\le C|\tau|^{-1-2s}
					&&(0<|\tau|\le1),\\
					C^{-1}e^{-p_s^*|\tau|}&\le k_s(\tau)\le Ce^{-p_s^*|\tau|}
					&&(|\tau|\ge1).
				\end{aligned}
			\end{equation}
			\item Let \(0<\sigma<2-2s\), \(0\le\lambda<{p_s^*}\), and let \(v\in C^{2s+\sigma}_{\rm loc}(\R)\) satisfy \(|v(\tau)|\le Ce^{\lambda|\tau|}\).  Then
			\[
				\begin{aligned}
					\Ps v(t)&\coloneqq\PV\int_{\R}(v(t)-v(\tau))k_s(t-\tau)\,d\tau+m_s v(t)\\
					&\;=\frac12\int_\R\bigl(2v(t)-v(t+\tau)-v(t-\tau)\bigr)k_s(\tau)\,d\tau+m_sv(t)
				\end{aligned}
			\]
			is well defined (the last integral converges absolutely), the function \(u(x)\coloneqq|x|^{-a}v(-a\log|x|)\) belongs to \(C^{2s+\sigma}_{\rm loc}(\R^N\setminus\{0\})\), the integral in \eqref{eq:frac-lap} converges for every \(x\neq0\), and
			\begin{equation}\label{eq:transform-general}
				(-\Delta)^su(x)=|x|^{-\frac{N+2s}{2}}\,\Ps v(t),
				\qquad t=-a\log|x|,\quad x\ne0 .
			\end{equation}
			\item Consequently, for \(v\) as in {\rm(ii)} with \(v>0\), the function \(u\) solves \((-\Delta)^su=u^p+u^{{p_s^*}+\eps}\) in \(\R^N\setminus\{0\}\) if and only if
			\begin{equation}\label{eq:transformed}
				\Ps v=e^{-bt}v^p+e^{\eps t}v^{{p_s^*}+\eps}\qquad\text{in }\R .
			\end{equation}
		\end{enumerate}
	\end{proposition}

	\begin{proof}
		(i) Evenness follows from \eqref{eq:K-inversion}. As \(\rho\to0\), \(K(\rho)=|\Sph^{N-1}|(1+O(\rho^2))\), because the linear term integrates to zero on the sphere; this gives \eqref{eq:kernel-infinity}. Near \(\rho=1\), write \(\rho=1+\zeta\) and use polar coordinates on \(\Sph^{N-1}\) around \(e_1\):
		\[
			K(1+\zeta)=|\Sph^{N-2}|\,I_{N,s}\,(1+\zeta)^{-\frac{N-1}2}|\zeta|^{-1-2s}+O\Bigl(|\zeta|^{1-2s}+\log\frac2{|\zeta|}\Bigr),
		\]
		where \(I_{N,s}=\frac12B\bigl(\frac{N-1}2,s+\frac12\bigr)\). Substitute \(\zeta=e^{\tau/a}-1\). Since \({p_s^*}-\frac{N-1}{2a}-\frac{1+2s}{2a}=0\), the term of order \(|\tau|^{-2s}\) cancels, giving \eqref{eq:kernel-expansion} with
		\begin{equation*}
			c_1=a^{2s}c_{1,s},\qquad
			c_{1,s}=\frac{4^ss\Gamma(s+1/2)}{\sqrt\pi\Gamma(1-s)}.
		\end{equation*}
		The remainder comes from replacing \(\sin^{N-2}\varphi\) and \(4(1+\zeta)\sin^2(\varphi/2)\) by \(\varphi^{N-2}\) and \((1+\zeta)\varphi^2\) on \([0,\pi/2]\), which produces an error bounded by a multiple of
		\[
			\int_0^{\pi/2}\frac{\varphi^N\,d\varphi}{(\zeta^2+\varphi^2)^{\frac{N+2s}2}}
			\le\int_0^{\pi/2}\frac{d\varphi}{(\zeta^2+\varphi^2)^{s}}
			\le C\left(|\zeta|^{1-2s}+\log\frac2{|\zeta|}\right).
		\]
		
		For \(s>\frac12\), the last inequality follows from the change of variables \(\varphi=|\zeta|\upsilon\), and for \(s\le\frac12\) from \((\zeta^2+\varphi^2)^{-s}\le\min\{|\zeta|^{-2s},\varphi^{-2s}\}\). Extending the rescaled main integral to infinity gives an \(O(1)\) error; for \(|\zeta|>\frac12\), the kernel is bounded. Thus \(\tilde k\in L^1\) and \((\cosh\tau-1)k_s\in L^1\), with \(0<m_s<\infty\).
		
		\medskip

		\noindent(ii) The second-difference integrand is \(O(|\tau|^{\sigma-1})\) near zero and \(O(e^{-({p_s^*}-\lambda)|\tau|})\) at infinity, so the integral converges absolutely. Identity \eqref{eq:transform-general} is the radial form of the conformal fractional Laplacian on \(\R\times\Sph^{N-1}\); see \cite{MR3694655,AoChanDelaTorreEtAl2019}. To fix the constants, take \(w\in C_c^\infty(\R)\) and 
		\[
			\tilde u(x)\coloneqq|x|^{-a}w(-a\log|x|).
		\] 
		The truncated operators 
		\[
			L_\delta u(x)\coloneqq c_{N,s}\int_{|x-y|>\delta}\frac{u(x)-u(y)}{|x-y|^{N+2s}}\,dy
		\] 
		converge to \((-\Delta)^su\) uniformly on \(\supp\tilde u\), by the second-difference formula. Symmetrization and the limit \(\delta\to0\) give 
		\[
			\int_{\R^N}(-\Delta)^su(x)\tilde u(x) \,dx=\frac{c_{N,s}}2\iint_{\R^{2N}}\frac{(u(x)-u(y))(\tilde u(x)-\tilde u(y))}{|x-y|^{N+2s}}\,dx\,dy.
		\]
		Dominated convergence applies: away from \(0\), the integrand is \(O(|x-y|^{1+\varkappa_0-N-2s})\) near the diagonal, where \(\varkappa_0\coloneqq\min\{1,2s+\sigma\}\) and \(1+\varkappa_0>2s\). The tail is controlled by \(|u(y)|(1+|y|)^{-N-2s}\in L^1\), since \(|v|\le Ce^{\lambda|\tau|}\) and \(\lambda<{p_s^*}\). No global \(H^s\) assumption is needed. In this absolutely convergent integral, set \(x=r\theta\), \(y=r\rho\omega\), \(t=-a\log r\), and \(\tau=t-a\log\rho\). The power of \(r\) is \(-1\), and
		\[
			\frac{c_{N,s}}a\rho^NK(\rho)=e^{t-\tau}k_s(t-\tau).
		\]
		Introduce the cylindrical cutoff \(|t-\tau|>\delta\) after this change of variables; it does not correspond to \(|x-y|>\delta\). All terms are now absolutely integrable. Near the diagonal, the quadratic part is \(O(|t-\tau|^{1+\varkappa_0})k_s(t-\tau)\). The zeroth-order terms are multiples of \(\int_\R vw\,dt\int_{|\zeta|>\delta}(e^{\pm\zeta}-1)k_s(\zeta)\,d\zeta\), whose odd parts vanish because \(k_s\) is even. Letting \(\delta\to0\) and using
		\[
			\lim_{\delta\to0}\int_{|\zeta|>\delta}(e^{\pm\zeta}-1)k_s(\zeta)\,d\zeta=m_s
		\] 
		gives
		\begin{equation}\label{eq:weak-EF}
			\begin{aligned}
				\int_{\R^N}(-\Delta)^su\,\tilde u\,dx
				&=\frac{|\Sph^{N-1}|}{a}\Bform(v,w),\\
				\Bform(v,w)
				&\coloneqq\frac12\iint_{\R^2}
				(v(t)-v(\tau))(w(t)-w(\tau))k_s(t-\tau)\,dt\,d\tau
				+m_s\int_\R vw\,dt,
			\end{aligned}
		\end{equation}
		and \(\Bform(v,w)=\int_\R(\Ps v)w\) by Fubini. Polar coordinates and \(N-a=\frac{N+2s}2\) give \eqref{eq:transform-general} in distributions, and pointwise by continuity.
		
		\medskip

		\noindent (iii) Use \(r^{a{p_s^*}}u^p=e^{-bt}v^p\) and \(r^{a{p_s^*}}u^{{p_s^*}+\eps}=e^{\eps t}v^{{p_s^*}+\eps}\).
	\end{proof}

	We solve \eqref{eq:transformed} with \(v^p\) and \(v^{{p_s^*}+\eps}\) replaced by \((v_+)^p\) and \((v_+)^{{p_s^*}+\eps}\). Lemma~\ref{lem:positive} gives positivity.

	\begin{remark}
		For \(v=1\), Proposition~\ref{prop:EF} gives the first identity. The value of \(m_s\) is the sharp fractional Hardy constant \cite[(2.1)--(2.2)]{FrankLiebSeiringer2008}; see also \cite[(1.2)]{MR3694655} for the cylindrical normalization:
		\[
			(-\Delta)^s|x|^{-a}=m_s|x|^{-a-2s},\qquad
			m_s=4^s\frac{\Gamma((N+2s)/4)^2}{\Gamma((N-2s)/4)^2}.
		\]
		Thus the zeroth-order coefficient of \(\Ps\) is the optimal constant in the fractional
		Hardy inequality
		\[
			m_s\int_{\R^N}\frac{|\chi(x)|^2}{|x|^{2s}}\,dx
			\le\int_{\R^N}|(-\Delta)^{s/2}\chi|^2\,dx,
			\qquad \chi\in C_c^\infty(\R^N).
		\]

	\end{remark}

	The energy uses \(\Bform\) from \eqref{eq:weak-EF}, including its factor \(1/2\). We have
	\[
		\Bform(v,w)=\int_\R(\Ps v)w\,dt
	\]
	for the functions in Proposition~\ref{prop:EF}(ii) and \(w\in C_c^\infty(\R)\). By density, 
	the identity also holds for \(v,w\in H_k^s(\R)\) whenever \(\Ps v\in L^2\); see Appendix~\ref{app:functional}.

	\medskip
	
	Throughout the paper, we use the Fourier convention
	\[
		\widehat f(k)=\int_\R e^{-ikt}f(t)\,dt,\qquad
		f(t)=\frac1{2\pi}\int_\R e^{ikt}\widehat f(k)\,dk.
	\]
	For \(w,f\in L^2(\R)+L^\infty(\R)\), the equation \(\Ps w=f\) in \(\mathcal D'(\R)\) means
	\[
		\int_\R w\Ps\varphi\,dt=\int_\R f\varphi\,dt
		\qquad(\varphi\in C_c^\infty(\R)).
	\]
	These integrals are finite since \(\Ps\varphi\in L^1\cap L^\infty\). By Lemma~\ref{lem:indicial}(ii), \(\Psym\) and all its derivatives are polynomially bounded. Thus \(\Ps\colon\mathcal S(\R)\to\mathcal S(\R)\) is continuous and extends by duality:
	\[
		\langle\Ps w,\varphi\rangle=\langle w,\Ps\varphi\rangle.
	\]
	The density of \(C_c^\infty(\R)\) in \(\mathcal S(\R)\) and the inclusion \(\mathcal S(\R)\hookrightarrow L^1(\R)\cap L^2(\R)\) extend the equation to test functions \(\varphi\in\mathcal S(\R)\). Therefore \(\Ps w=f\) in \(\mathcal D'(\R)\) if and only if it holds in \(\mathcal S'(\R)\).

	\begin{lemma}[Indicial roots and the Green kernel]\label{lem:indicial}
		The operator \(\Ps\) has the following properties.
		\begin{enumerate}[label=\textup{(\roman*)}]
			\item For \(|\lambda|<{p_s^*}\),
			\(\Ps e^{\lambda t}=\Theta(\lambda)e^{\lambda t},\)
			where
			\[
				\Theta(\lambda)\coloneqq\int_\R
				(\cosh\tau-\cosh(\lambda\tau))k_s(\tau)\,d\tau.
			\]
			The function \(\Theta\) is even and strictly decreasing on \([0,{p_s^*})\), with \(\Theta(0)=m_s\) and \(\Theta(\pm1)=0\).
			\item The Fourier symbol is
			\[
				\Psym(k)=m_s+\int_\R(1-\cos(k\tau))k_s(\tau)\,d\tau.
			\]
			It is smooth and even on \(\R\), and
			\[
				c(1+|k|^{2s})\le\Psym(k)\le C(1+|k|^{2s}),\qquad
				|\partial_k^j\Psym(k)|\le C_j(1+|k|^{2s}).
			\]
			For \(\varphi\in\mathcal S(\R)\), \(\widehat{\Ps\varphi}=\Psym\widehat\varphi\).
			The symbol extends analytically to \(|\Ima k|<{p_s^*}\), where
			\begin{equation}\label{eq:Re-symbol}
				\Rea\Psym(k_1+i\lambda)
				=\Theta(\lambda)+\int_\R
				(1-\cos(k_1\tau))\cosh(\lambda\tau)k_s(\tau)\,d\tau.
			\end{equation}
			Its only zeros in \(|\Ima k|\le1\) are \(\pm i\). They are simple, with
			\[
				\Psym'(\pm i)=\pm i\beta_s,\qquad
				\beta_s\coloneqq\int_\R\tau\sinh\tau\,k_s(\tau)\,d\tau=-\Theta'(1)>0.
			\]
			\item The Green kernel \(G_s\coloneqq\Fou^{-1}(1/\Psym)\) is an even, strictly positive function in \(L^1(\R)\), continuous on \(\R\setminus\{0\}\). Set
			\[
				c^{\mathrm R}_{N,s}\coloneqq\frac{\Gamma(a)}{4^s\pi^{N/2}\Gamma(s)}.
			\]
			Then
			\begin{equation}\label{eq:Green-Riesz}
				G_s(z)=\frac{c^{\mathrm R}_{N,s}}a e^z
				\int_{\Sph^{N-1}}\frac{d\omega}
				{|e_1-e^{z/a}\omega|^{N-2s}},
			\end{equation}
			for \(z\ne0\). In particular,
			\[
				\begin{gathered}
					G_s(z)=\beta_s^{-1}e^{-|z|}
					\bigl(1+O(e^{-2|z|/a})\bigr)\quad(|z|\to\infty),\\
					G_s(z)\le Ce^{-|z|}\quad(|z|\ge1),\qquad C^{-1}\varpi_s(z)\le G_s(z)\le C\varpi_s(z)\quad(0<|z|\le1).
				\end{gathered}
			\]
			The coefficient is
			\[
				\beta_s^{-1}=\frac{c^{\mathrm R}_{N,s}|\Sph^{N-1}|}{a}
				=\frac{2\Gamma(a)}{a4^s\Gamma(s)\Gamma(N/2)}.
			\]
			Moreover, \(G_s\in C^1(\R\setminus\{0\})\), with
				\[
					|G_s'(z)|\le C|z|^{2s-2}\quad(0<|z|\le1),\qquad|G_s'(z)|\le Ce^{-|z|}\quad(|z|\ge1),
				\]
				\(G_s\in L^{\varsigma}(\R)\) for every \(1\le\varsigma<\varsigma_s\), and \(G_s*f\) is bounded and uniformly continuous for every \(f\in L^\infty(\R)\). If \(s>\frac12\), then \(G_s\) is bounded and continuous on \(\R\), while if \(0<s\le\frac12\), then \(G_s\) is unbounded near \(0\).
			\item If \(w,f\in L^2(\R)\) and \(\Ps w=f\) in \(\mathcal D'(\R)\), then \(w=G_s*f\). If \(w\in L^\infty(\R)\) and \(\Ps w=0\) in \(\mathcal D'(\R)\), then \(w=0\).
		\end{enumerate}
	\end{lemma}

	\begin{proof}
		(i) For \(|\lambda|<{p_s^*}\), the odd part cancels in the principal value, giving \(\Ps e^{\lambda t}=\Theta(\lambda)e^{\lambda t}\). Directly from the definitions,
		\[
			\Theta(0)=m_s,\qquad \Theta(\pm1)=0,\qquad
			\Theta'(\lambda)=-\int_\R\tau\sinh(\lambda\tau)k_s(\tau)\,d\tau<0\quad(\lambda>0).
		\]
		The modes \(e^{-t}\) and \(e^t\) correspond to \(1\) and \(|x|^{2s-N}\).
		
		\medskip

		\noindent (ii) Take the Fourier transform of the second-difference formula. The bounds follow from \eqref{eq:kernel-bounds}, splitting the integral at \(|\tau|=|k|^{-1}\); for the first derivative, we use \(|\tau\sin(k\tau)|\le|\tau|\min\{1,|k||\tau|\}\). In the strip, \eqref{eq:Re-symbol} and \(\Theta>0\) on \((-1,1)\) show that the only zeros with \(|\Ima k|\le1\) are \(\pm i\), and \(\Psym'(\pm i)=\pm i\beta_s\).
		
		\medskip
		
		\noindent (iii) The multipliers \(\Psym\) and \(1/\Psym\) are smooth with polynomially bounded derivatives. Hence, for \(h\in\mathcal S(\R)\),
		\[
			G_s*h\coloneqq\Fou^{-1}(\widehat h/\Psym)\in\mathcal S(\R),
			\qquad \Ps(G_s*h)=h.
		\]
		Since \(\Psym\ge m_s>0\) on \(\R\), multiplication by \(1/\Psym\) is valid in \(\mathcal S'(\R)\); thus \(\Ps w=0\) implies \(w=0\) there. For \(h\in C_c^\infty(\R)\), set \(F(x)\coloneqq|x|^{-a-2s}h(-a\log|x|)\) and let
		\[
			u\coloneqq c^{\mathrm R}_{N,s}|\cdot|^{-(N-2s)}*F
		\]
		be its Riesz 
		potential \cite[Ch.~V, \S1]{Stein1970}. Since \(F\) is smooth and supported in an annulus, \((-\Delta)^su=F\), \(u\) is smooth, and \(v=r^au=O(e^{-|t|})\) at both ends. Proposition~\ref{prop:EF}(ii) with \(\lambda=0\) gives \(\Ps v=h\) pointwise. For \(\delta>0\), define the cylindrical truncation
		\[
			\Ps^\delta\varphi(t)\coloneqq
			\int_{|z|>\delta}\bigl(\varphi(t)-\varphi(t-z)\bigr)k_s(z)\,dz+m_s\varphi(t).
		\]
		For bounded \(w\) and \(\varphi\in C_c^\infty(\R)\), Fubini's theorem (since \(\int_{|z|>\delta}k_s(z)\,dz<\infty\)) and the evenness of \(k_s\) give
		\(\int(\Ps^\delta w)\varphi=\int w\Ps^\delta\varphi\).
		The second-difference formula and \eqref{eq:kernel-bounds} give
		\(\|\Ps^\delta\varphi-\Ps\varphi\|_{L^1}\le C\delta^{2-2s}\|\varphi''\|_{L^1}\) for \(0<\delta<1\).
		If also \(w\in C^{2s+\sigma}_{\rm loc}(\R)\), with \(0<\sigma<2-2s\), then \(\Ps^\delta w\to\Ps w\) locally uniformly. Passing to the limit yields
		\(\int(\Ps w)\varphi=\int w\Ps\varphi\).
		Thus \(\Ps v=h\) in \(\mathcal D'(\R)\), with uniqueness among bounded solutions as above. After the substitution \(y=r\rho\omega\), \(z=a\log\rho\), the powers of \(r\) cancel and \(v=\widetilde G_s*h\), where \(\widetilde G_s\) is the right-hand side of \eqref{eq:Green-Riesz}. For \(\frac12\le\rho\le2\), 
		\[
			|e_1-\rho\omega|^2=(1-\rho)^2+\rho|e_1-\omega|^2\ge\frac12\bigl((1-\rho)^2+|e_1-\omega|^2\bigr),
		\]
		and \(|e_1-\omega|\ge2\varphi/\pi\) if 
		\(\omega_1=\cos\varphi\). Hence
		\[
			\int_{\Sph^{N-1}}\frac{d\omega}{|e_1-\rho\omega|^{N-2s}}
			\le C\int_0^\pi(|1-\rho|+\varphi)^{2s-2}\,d\varphi
			\le C\varpi_s(1-\rho).
		\]
		Thus \(\widetilde G_s(z)\le C\varpi_s(z)\) for \(0<|z|\le1\), since \(|1-e^{z/a}|\ge c|z|\). For the lower bound, assume
		\(0<|1-\rho|\le\frac14\). Then \(|e_1-\rho\omega|^2\le2\bigl((1-\rho)^2+|e_1-\omega|^2\bigr)\) and
		\(|e_1-\omega|\le\varphi\), so
		\[
			\int_{\Sph^{N-1}}\frac{d\omega}{|e_1-\rho\omega|^{N-2s}}
			\ge c\int_{|1-\rho|}^{\pi/2}\frac{\varphi^{N-2}\,d\varphi}{(|1-\rho|+\varphi)^{N-2s}}
			\ge c\int_{|1-\rho|}^{\pi/2}\varphi^{2s-2}\,d\varphi
			\ge c\,\varpi_s(1-\rho).
		\]
		Since \(|1-e^{z/a}|\le C|z|\) for \(|z|\le1\) and \(\varpi_s\) is nonincreasing in \(|z|\), this gives \(\widetilde G_s(z)\ge c\varpi_s(z)\) for \(0<|z|\le z_0\), with some \(z_0>0\). Continuity and positivity on \(z_0\le|z|\le1\) extend the bound to \(0<|z|\le1\). In particular, \(\widetilde G_s\) is unbounded near zero if \(s\le\frac12\). For the derivative, set 
		\[
			I(\rho)\coloneqq\int_{\Sph^{N-1}}|e_1-\rho\omega|^{2s-N}\,d\omega,
		\] 
		so that 
		\[
			\widetilde G_s'(z)=\frac{c^{\mathrm R}_{N,s}}ae^z\bigl(I(\rho)+\frac\rho aI'(\rho)\bigr)\quad  \text{with } \rho=e^{z/a}.
		\] 
		Since 
		\(|\rho-\omega_1|\le|e_1-\rho\omega|\), we have 
		\[
			|I'(\rho)|\le(N-2s)\int_{\Sph^{N-1}}|e_1-\rho\omega|^{2s-N-1}\,d\omega.
		\] 
		The preceding estimate gives, for \(\frac12\le\rho\le2\),
		\[
			|I'(\rho)|\le C\int_0^\pi(|1-\rho|+\varphi)^{2s-3}\,d\varphi\le C|1-\rho|^{2s-2}.
		\]
		Together with \(I(\rho)\le C\varpi_s(1-\rho)\le C|1-\rho|^{2s-2}\), this yields \(|\widetilde G_s'(z)|\le C|z|^{2s-2}\) for \(0<|z|\le1\). Since \(I\) and \(I'\) are bounded on compact subsets of \([0,\infty)\setminus\{1\}\), we also have \(|\widetilde G_s'(z)|\le Ce^z\) for \(z\le-1\); evenness gives the bound for \(z\ge1\). Thus \(G_s=\widetilde G_s\) in \(\mathcal D'(\R)\) and as functions. This kernel is positive, even by inversion, and continuous away from zero. As \(z\to-\infty\),
		\[
			\widetilde G_s(z)=C_se^{z}(1+O(e^{2z/a}))
		\]
		with \(C_s=c^{\mathrm R}_{N,s}|\Sph^{N-1}|/a\). Finally, 
		\[
			\int e^{\lambda t}G_s=1/\Theta(\lambda)\quad\text{ for }0<\lambda<1
		\] 
		by analytic continuation, and letting \(\lambda\uparrow1\) gives \(C_s=\lim(1-\lambda)/\Theta(\lambda)=\beta_s^{-1}\). If \(s>\frac12\), the bound on \(G_s'\) is integrable near zero. Thus the finite one-sided limits of \(G_s\) coincide by evenness and define its continuous extension at zero. The remaining assertions follow from \(G_s\in L^1\cap L^{\varsigma}\) and the continuity of translations in \(L^1\) and \(L^{\varsigma}\).

		\medskip
		
		\noindent(iv) Since \(\widehat{G_s}=1/\Psym\), Plancherel gives the first assertion. For the second, the equivalence of the \(\mathcal D'\) and \(\mathcal S'\) formulations gives \(\Ps w=0\) in \(\mathcal S'(\R)\). Part (iii) then gives \(w=0\).
	\end{proof}

	For \(0\le\beta<1<\eta\) and \(t,\xi\in\R\),
	\begin{equation}\label{eq:conv-elementary}
		\begin{aligned}
			&\int_\R e^{-|t-r|}e^{-\beta|r-\xi|}\,dr
			\le C_\beta e^{-\beta|t-\xi|},\\
			&\int_\R e^{-|t-r|}e^{-|r-\xi|}\,dr
			=(1+|t-\xi|)e^{-|t-\xi|},\\
			&\int_\R e^{-|t-r|}e^{-\eta|r-\xi|}\,dr
			\le C_\eta e^{-|t-\xi|}.
		\end{aligned}
	\end{equation}
	These estimates follow by splitting the integral at \(r=\xi\) and \(r=t\); for \(\lambda\ne1\), 
	\[
		\int_\R e^{-|t-r|}e^{-\lambda|r|}\,dr=\frac{2}{1-\lambda^2}(e^{-\lambda|t|}-\lambda e^{-|t|}).
	\]

	For \(s\le\frac12\), the bound \(G_s(z)\le Ce^{-|z|}\) holds only for \(|z|\ge1\). To estimate its convolution with the majorants used below, call a measurable function \(g\colon\R\to[0,\infty)\) \emph{\(M\)-slowly varying} if
	\begin{equation}\label{eq:slowly-varying}
		g(r)\le Mg(t)\qquad\text{whenever }|r-t|\le1.
	\end{equation}

	\begin{lemma}[Comparison with the exponential kernel]\label{lem:comparison}
		There is \(C>0\), depending only on \(N\) and \(s\), such that every \(M\)-slowly varying function \(g\) satisfies
		\[
			(G_s*g)(t)\le C(1+M^2)\int_\R e^{-|t-r|}g(r)\,dr
			\qquad(t\in\R).
		\]
	\end{lemma}

	\begin{proof}
		By Lemma~\ref{lem:indicial}(iii),
		\[
			\int_{|t-r|\ge1}G_s(t-r)g(r)\,dr\le C\int_\R e^{-|t-r|}g(r)\,dr.
		\]
		By \eqref{eq:slowly-varying},
		\[\begin{gathered}
			\int_{|t-r|<1}G_s(t-r)g(r)\,dr\le Mg(t)\|G_s\|_{L^1},\\
			\text{ and }
			\int_\R e^{-|t-r|}g(r)\,dr\ge e^{-1}\int_{|t-r|\le1}g(r)\,dr\ge\frac2{eM}g(t).\end{gathered}
		\]
		The result follows by addition.
	\end{proof}

	All majorants below are slowly varying, with constants depending only on \(n\), \(\beta\), \(\gamma\), \(\eta\), \(\theta_0\) and \({p_s^*}\), independently of the centers. The functions \(e^{-\beta|t-\xi|}\) and \((1+|t-\xi|)e^{-|t-\xi|}\) have constants \(M\le2e^\beta\) and \(M\le2e\); finite sums and products preserve this property. This includes \(\Wgt_\beta\), \(Q_\pm\), and the majorants defining \(X\) and \(Y\) in Subsection~\ref{subsec:tower}: the pieces \(\Wgt_\gamma\) on \([\xi_1,\xi_n]\) and \(Q_\pm\), \(e^{-\eta(\xi_1-t)}\), or \(e^{-\eta(t-\xi_n)}\) on the exterior half-lines are comparable at the endpoints up to the factor \(n\). Thus, if \(|f|\le g\) with \(g\) slowly varying, then
	\[
		|G_s*f|\le C\,e^{-|\cdot|}*g,
	\]
	which is estimated by \eqref{eq:conv-elementary}. Each bootstrap uses a fixed number of steps, so the constants remain independent of \(\eps\) and the centers.

\subsection{The critical profile and nondegeneracy}

	Set
	\[
		W(x)\coloneqq W_{1,0}(x)=\alpha_{N,s}(1+|x|^2)^{-a},
	\]
	and define
	\[
		U(t)\coloneqq e^{-t}W(e^{-t/a}).
	\]
	Since \(e^{-t}(1+e^{-2t/a})^{-a}=(e^{t/a}+e^{-t/a})^{-a}\),
	\begin{align*}
		U(t)
		&=\alpha_{N,s}e^{-t}
		\left(1+e^{-\frac{2}{a}t}\right)^{-a}
		=\alpha_{N,s}\left(2\cosh\frac{t}{a}\right)^{-a}.
	\end{align*}
	Thus \(U\) is even, \(U'\) is odd, and
	\[
		U(t)=\alpha_{N,s}e^{-|t|}
		\bigl(1+O(e^{-2|t|/a})\bigr)\qquad(|t|\to\infty).
	\]
	Also,
	\[
			U'(t)=-\tanh(t/a)U(t),\quad
			U''(t)=\left[\tanh^2(t/a)-a^{-1}\operatorname{sech}^2(t/a)\right]U(t).
	\]
	The coefficients in these identities have bounded derivatives of every order. Hence
	\begin{equation}\label{eq:U-tail}
		ce^{-|t|}\le U(t)\le Ce^{-|t|},\qquad
		|U^{(j)}(t)|\le C_je^{-|t|}\quad(j\ge1).
	\end{equation}
	Since \(W\) solves the critical equation, Proposition~\ref{prop:EF} yields
	\begin{equation}\label{eq:Ueq}
		\Ps U=U^{p_s^*}\qquad\text{in }\R.
	\end{equation}

	Translations of \(U\) correspond to dilations in the original variables. If
	\[
		U_\xi(t)\coloneqq U(t-\xi),\qquad \mu=e^{-\xi/a},
	\]
	then
	\begin{equation}\label{eq:translation-scale}
		e^tU(t-\xi)=W_{\mu,0}(r),\qquad r=e^{-t/a}.
	\end{equation}
	Differentiating \eqref{eq:Ueq} shows that \(U'\) solves
	\begin{equation}\label{eq:lin-critical}
		\Ps\psi={p_s^*} U^{{p_s^*}-1}\psi\qquad\text{in }\R ,
	\end{equation}
	while the scaling derivative satisfies
	\begin{equation}\label{eq:scaling-mode}
		x\cdot\nabla W+aW=-a\,e^{t}U'(t),\qquad t=-a\log|x| .
	\end{equation}

	The following cutoff argument removes the singularity at the origin.

	\begin{lemma}[Removable singularity at the origin]\label{lem:removable}
		Let \(\phi\in L^\infty(\R^N)\cap C^{2s+\sigma}_{\rm loc}(\R^N\setminus\{0\})\) for some \(0<\sigma<2-2s\), and let \(g\in L^\infty_{\rm loc}(\R^N)\).  If \((-\Delta)^s\phi=g\) pointwise in \(\R^N\setminus\{0\}\), then
		\[
			\int_{\R^N}\phi\,(-\Delta)^s\varphi\,dx=\int_{\R^N}g\,\varphi\,dx
			\qquad\text{for all }\varphi\in C_c^\infty(\R^N).
		\]
		That is, \((-\Delta)^s\phi=g\) in \(\mathcal D'(\R^N)\).
	\end{lemma}

	\begin{proof}
		For \(\varphi\in C_c^\infty(\R^N\setminus\{0\})\), truncating the kernel gives 
		\[
			\int_{\R^N}\phi(-\Delta)^s\varphi=\int_{\R^N}(-\Delta)^s\phi\,\varphi.
		\]
		This uses boundedness and local \(C^{2s+\sigma}\) regularity, without radial symmetry. Thus \(T\coloneqq(-\Delta)^s\phi-g\) is supported at zero. For \(0<\lambda\le1\), take a cutoff \(\zeta_\lambda=\zeta(\cdot/\lambda)\) and write
		\[
			\varphi\zeta_\lambda=\varphi(0)\zeta_\lambda+\chi_\lambda,
			\qquad \langle T,\varphi\rangle=\langle T,\varphi\zeta_\lambda\rangle.
		\]
		Scaling gives
		\[
			\|(-\Delta)^s\zeta_\lambda\|_{L^1}=C\lambda^{N-2s},
			\qquad \|(-\Delta)^s\chi_\lambda\|_{L^1}\le C\lambda^{N+1-2s}.
		\]
		Since \(N>2s\),
		\[
			|\langle T,\varphi\rangle|\le C\bigl(\lambda^{N-2s}+\lambda^{N+1-2s}+\lambda^N\bigr)\longrightarrow0
			\qquad(\lambda\to0).
		\]
	\end{proof}

	\begin{lemma}[Radial nondegeneracy]\label{lem:nondeg}
		Let \(\psi\in L^\infty(\R)\) solve \eqref{eq:lin-critical} in \(\mathcal D'(\R)\).  Then \(\psi\in C^{2s+\sigma}_{\rm loc}(\R)\) for some \(0<\sigma<2-2s\), \(|\psi(t)|\le Ce^{-|t|}\) for all \(t\in\R\), and \(\psi=cU'\) for some \(c\in\R\).
	\end{lemma}

	\begin{proof}
		\noindent\emph{Step 1.} Decay and regularity. Set
		\[
			f\coloneqq{p_s^*}U^{{p_s^*}-1}\psi.
		\]
		We have \(f=O(e^{-({p_s^*}-1)|t|})\), hence \(f\in L^1\cap L^2\cap L^\infty\). By Lemma~\ref{lem:indicial}(iii), \(G_s*f\) is bounded, and \(\widehat{G_s}=1/\Psym\) gives \(\Ps(G_s*f)=f\). Thus \(\psi-G_s*f\) is a bounded homogeneous solution, and Lemma~\ref{lem:indicial}(iv) yields
		\[
			\psi=G_s*f.
		\]
		By Lemma~\ref{lem:comparison} and \eqref{eq:conv-elementary}, each iteration increases the decay exponent by \({p_s^*}-1\), up to \(1\). If the exponent \(1\) is reached exactly, the bound \((1+|t|)e^{-|t|}\) improves to \(Ce^{-|t|}\) after one more step. Since \({p_s^*}U^{{p_s^*}-1}\) is smooth, a finite bootstrap with Proposition~\ref{prop:regularity} gives the local regularity.

		\medskip

		\noindent\emph{Step 2.} The original variables. Set \(\phi(x)\coloneqq|x|^{-a}\psi(-a\log|x|)\). By Step~1, \(|\phi(x)|\le C(1+|x|)^{2s-N}\). Proposition~\ref{prop:EF} and Lemma~\ref{lem:removable} give
		\[
			(-\Delta)^s\phi={p_s^*}W^{{p_s^*}-1}\phi\qquad\text{in }\mathcal D'(\R^N).
		\]
		By Step~1 and \eqref{eq:Riesz-cylinder}, 
		\[
			\phi=c^{\mathrm R}_{N,s}|\cdot|^{2s-N}*({p_s^*}W^{{p_s^*}-1}\phi).
		\] 
		Since \(|{p_s^*}W^{{p_s^*}-1}\phi|\le C(1+|x|)^{-N-2s}\), the source \({p_s^*}W^{{p_s^*}-1}\phi\) belongs to \(L^1(\R^N)\cap L^\infty(\R^N)\). Choose \(\beta\in(0,\min\{1,2s,2-2s\})\) with \(2s+\beta\ne1\). Lemma~\ref{lem:Riesz} gives \(\phi\in C^{0,\beta}(\R^N)\). Since \(W^{{p_s^*}-1}\) is bounded and Lipschitz,
		\[
			{p_s^*}W^{{p_s^*}-1}\phi\in C^{0,\beta}(\R^N),
			\qquad \phi\in C^{2s+\beta}(\R^N)
		\]
		by \eqref{eq:Silvestre-holder}. Thus \(\phi\) is a bounded classical solution of the linearized equation on \(\R^N\).

		\medskip

		\noindent\emph{Step 3.} Conclusion. By \cite[Theorem~1.1]{DavilaDelPinoSire2013}, valid for every \(0<s<1\) and \(N>2s\),
		\[
			\phi=\sum_jc_j\partial_jW+c_0(x\cdot\nabla W+aW).
		\]
		Radial symmetry forces \(c_j=0\), and \eqref{eq:scaling-mode} gives \(\psi=-ac_0U'\).
	\end{proof}

\subsection{Tower parameters}\label{subsec:tower}

	Fix \(n\ge1\).  For \(\xi=(\xi_1,\ldots,\xi_n)\) with
	\[
		0<\xi_1<\xi_2<\cdots<\xi_n,
	\]
	set
	\[
		U_i(t)\coloneqq U(t-\xi_i),\qquad Z_i(t)\coloneqq U_i'(t),
		\qquad V(t)\coloneqq\sum_{i=1}^nU_i(t).
	\]
	Parametrize the centers by \(\Lambda=(\Lambda_1,\ldots,\Lambda_n)\), as in the local mixed-power problem \cite[(14)]{MR2381679}:
	\begin{equation}\label{eq:xigap}
		\xi_1=-\frac1b\log\eps-\log\Lambda_1,\quad
		\xi_i-\xi_{i-1}=-\log\eps-\log\Lambda_i,
		\qquad 2\le i\le n.
	\end{equation}
	Fix a small \(\delta\in(0,1)\) and restrict the parameters to
	\[
		\mathcal K_\delta\coloneqq\{\Lambda\in\R^n:\delta\le\Lambda_i\le\delta^{-1}\ \text{for all }i\}.
	\]
	Uniformly on \(\mathcal K_\delta\),
	\begin{equation}\label{eq:separation}
		\xi_1\to\infty,\qquad
		\min_{2\le i\le n}(\xi_i-\xi_{i-1})\to\infty\quad(n\ge2),
		\qquad
		\eps\xi_n\to0\qquad(\eps\downarrow0).
	\end{equation}
	Iteration of \eqref{eq:xigap} gives
	\begin{equation}\label{eq:xi-explicit}
		\xi_i
		=-\left(\frac1b+i-1\right)\log\eps
		-\sum_{j=1}^i\log\Lambda_j,
		\qquad 1\le i\le n.
	\end{equation}
	Hence, uniformly for \(\Lambda\in\mathcal K_\delta\),
	\begin{align}
		&e^{-b\xi_1}=\eps\Lambda_1^b,\notag\\
		&e^{-(\xi_i-\xi_{i-1})}=\eps\Lambda_i,
		\qquad 2\le i\le n,\label{eq:int-balance}\\
		&\eps\xi_i=O(\eps|\log\eps|)=o(1).\notag
	\end{align}
	The first two identities fix the subcritical and interaction terms at order \(\eps\). The change of parameters is invertible, with triangular Jacobian
	\[
		\frac{\partial\xi_i}{\partial\Lambda_j}
		=\begin{cases}-\Lambda_j^{-1},&j\le i,\\0,&j>i.\end{cases}
	\]

	Fix exponents
	\begin{equation}\label{eq:gamma-eta}
		\begin{gathered}
			0<\gamma<\min\{1,{p_s^*}-1\},
			\qquad
			1<\eta<2p-{p_s^*},\\
			\text{and set}\qquad
			\theta_0\coloneqq2p-{p_s^*}-1>0 .
		\end{gathered}
	\end{equation}
	By \eqref{eq:2pP}, the interval for \(\eta\) is nonempty. Since \(p<{p_s^*}\),
	\[
		2p-{p_s^*}<{p_s^*},\qquad
		\theta_0=\min\{{p_s^*}-1,2p-{p_s^*}-1\},\qquad \eta<{p_s^*}.
	\]
	These choices ensure that:
	\begin{itemize}
		\item the decay exponent \(\gamma\) lies below the positive indicial root, and the interaction remainder is smaller than \(\eps^{1/2}\);
		\item the exterior sources decay faster than the Green kernel and at a rate compatible with the subcritical term at minus infinity;
		\item the linearized potential decays away from the centers at rate \(\theta_0>0\).
	\end{itemize}
	In particular, \(\eta<1+\theta_0\). In the proof of Theorem~\ref{thm:main}, we take \(\eta\coloneqq\frac12\bigl(1+2p-{p_s^*}\bigr)\), so that \(\eta\) depends only on \(N\), \(s\) and \(p\). For measurable functions, we define
	\[
		\Wgt_\beta(t)\coloneqq\sum_{i=1}^n e^{-\beta|t-\xi_i|},
		\qquad
		\|f\|_\beta\coloneqq\operatorname*{ess\,sup}_{t\in\R}
		\Wgt_\beta(t)^{-1}|f(t)|,\qquad \beta>0 .
	\]
	In particular, \(\|\cdot\|_1\) denotes this norm with \(\beta=1\), not the \(L^1\) norm.

	The weights satisfy
	\begin{align*}
		\Wgt_\beta&\le n,
		&\Wgt_{\beta_1}&\le\Wgt_{\beta_2}\quad(\beta_1\ge\beta_2),\\
		\Wgt_{\beta_1}\Wgt_{\beta_2}&\le n^2\Wgt_{\beta_1+\beta_2},
		&\Wgt_\beta^{\,\nu}&\le n^{\nu-1}\Wgt_\beta\quad(\nu\ge1).
	\end{align*}
	These estimates follow from \(e^{-\beta\dist(t,\{\xi_i\})}\le\Wgt_\beta(t)\le ne^{-\beta\dist(t,\{\xi_i\})}\).
	By \eqref{eq:U-tail}, uniformly under \eqref{eq:separation},
	\[
		c\Wgt_1(t)\le V(t)\le C\Wgt_1(t).
	\]

	Lemma~\ref{lem:comparison} and \eqref{eq:conv-elementary} give \(|(G_s*Z_i)(t)|\le C(1+|t-\xi_i|)e^{-|t-\xi_i|}\). On the exterior half-lines, set
	\[
		Q_-(t)\coloneqq(1+\xi_1-t)e^{t-\xi_1},\qquad t\le\xi_1,\quad
		Q_+(t)\coloneqq(1+t-\xi_n)e^{-(t-\xi_n)},\qquad t\ge\xi_n,
	\]
	and define
	\begin{align*}
		\|f\|_{\Lout}
		&\coloneqq\max\bigg\{
		\sup_{t\le\xi_1}\frac{|f(t)|}{Q_-(t)},
		\sup_{t\ge\xi_n}\frac{|f(t)|}{Q_+(t)}
		\bigg\},\\
		\|h\|_{\eta,\Lout}
		&\coloneqq\max\bigg\{
		\operatorname*{ess\,sup}_{t\le\xi_1}e^{\eta(\xi_1-t)}|h(t)|,
		\operatorname*{ess\,sup}_{t\ge\xi_n}e^{\eta(t-\xi_n)}|h(t)|
		\bigg\}.
	\end{align*}
	Define the weighted spaces and their norms by
	\begin{equation*}
		\begin{aligned}
			&\|\phi\|_X\coloneqq\|\phi\|_\gamma+\|\phi\|_{\Lout},\qquad
			&&X\coloneqq\{\phi\in C(\R)\colon\|\phi\|_X<\infty\},	\\
			&\|h\|_Y\coloneqq\|h\|_\gamma+\|h\|_{\eta,\Lout},\qquad
			&&Y\coloneqq\{h\in L^\infty_{\rm loc}(\R)\colon\|h\|_Y<\infty\}.
		\end{aligned}	
	\end{equation*}
	Functions in \(Y\) are identified up to equality almost everywhere. Both spaces are
	Banach spaces. The weights depend on \(\xi\); we indicate this dependence in the
	subscript when comparing configurations.
	For \(y\ge0\) we write \(t_-(y)\coloneqq\xi_1-y\) and \(t_+(y)\coloneqq\xi_n+y\). The \emph{cell} of \(\xi_i\) is \(\Omega_i\coloneqq\{t\in\R:|t-\xi_i|=\dist(t,\{\xi_1,\ldots,\xi_n\})\}\); the cells cover \(\R\), and \(\Omega_1\supset(-\infty,\xi_1]\), \(\Omega_n\supset[\xi_n,\infty)\).

	\begin{lemma}[Moving the centers]\label{lem:moving-norms}
		Let \(\xi,\xi'\in\R^n\) be ordered configurations with \(|\xi-\xi'|\le1\).
		The corresponding \(X\), \(Y\), and \(\gamma\) norms are uniformly equivalent:
		\[
			\begin{aligned}
				C^{-1}\|f\|_{X,\xi}&\le\|f\|_{X,\xi'}\le C\|f\|_{X,\xi},\\
				C^{-1}\|h\|_{Y,\xi}&\le\|h\|_{Y,\xi'}\le C\|h\|_{Y,\xi},\\
				e^{-\gamma}\|h\|_{\gamma,\xi}
				&\le\|h\|_{\gamma,\xi'}\le e^\gamma\|h\|_{\gamma,\xi}.
			\end{aligned}
		\]
		The constant \(C\) depends only on \(n\), \(\gamma\) and \(\eta\), independently of the distances
		between consecutive centers. The same assertion holds for the weighted norms
		\(\|\cdot\|_{1,\xi}\) and \(\|\cdot\|_{1,\xi'}\).
	\end{lemma}

	\begin{proof}
		We have
		\[
			e^{-\beta}\Wgt_{\beta,\xi}\le\Wgt_{\beta,\xi'}\le e^\beta\Wgt_{\beta,\xi}.
		\]
		On a common exterior half-line, the endpoint distances satisfy \(|y-y'|\le1\), so
		\[
			\frac1{2e}\le\frac{(1+y')e^{-y'}}{(1+y)e^{-y}}\le2e,
			\qquad e^{-\eta}\le e^{\eta(y'-y)}\le e^\eta.
		\]
		The strip belonging to only one exterior half-line has length at most one. There the exterior weights are comparable to constants and \(\Wgt_\gamma\ge e^{-\gamma}\).
	\end{proof}

	Thus we differentiate with respect to the centers in fixed Banach spaces.

\section{The projected linear problem}\label{sec:linear}

	Fix \(n\ge1\), \(\delta\in(0,1)\) and \(\gamma,\eta\) as in \eqref{eq:gamma-eta}. For \(\Lambda\in\mathcal K_\delta\), let \(\xi\) be given by \eqref{eq:xi-explicit}. Solutions lie in \(X\cap H^s_k(\R)\), and right-hand sides in \(Y\). Constants may depend on the fixed parameters and are uniform for \(\Lambda\in\mathcal K_\delta\) and \(0<\eps<\eps_0\).

	The linearized operator at \(V\) is
	\[
		\mathfrak L_{\eps,\xi}\phi\coloneqq\Ps\phi-F_{\eps,\xi}\phi,
		\qquad
		F_{\eps,\xi}\coloneqq({p_s^*}+\eps)e^{\eps t}V^{{p_s^*}+\eps-1}
		+p e^{-bt}V^{p-1}.
	\]
	The projected problem is
	\begin{equation}\label{eq:linear-projected}
		\begin{cases}
			\mathfrak L_{\eps,\xi}\phi=h+\displaystyle\sum_{i=1}^n c_iZ_i
			&\text{in }\R,\\[2mm]
			\displaystyle\int_\R\phi Z_i=0,&i=1,\ldots,n.
		\end{cases}
	\end{equation}
	We use the energy space
	\begin{align*}
		H_k^s(\R)&\coloneqq\left\{w\in L^2(\R):
		\iint_{\R^2}|w(t)-w(r)|^2k_s(t-r)\,dt\,dr<\infty\right\},\\
		\|w\|_{H_k^s}^2&\coloneqq\Bform(w,w).
	\end{align*}
	A weak solution satisfies
	\[
		\Bform(\phi,\psi)-\int_\R F_{\eps,\xi}\phi\psi
		=\int_\R h\psi+\sum_{i=1}^n c_i\int_\R Z_i\psi,
		\qquad \psi\in H_k^s(\R).
	\]
	By Proposition~\ref{prop:H-equivalence}, \(H_k^s(\R)=H^s(\R)\) with equivalent norms, and the weak equation also holds in \(\mathcal D'(\R)\).

\subsection{The approximate kernel and the potential}

	Differentiating the equation for \(U_i\) gives
	\[
		\Ps Z_i={p_s^*} U_i^{{p_s^*}-1}Z_i.
	\]
	Thus
	\[
		\mathfrak L_{\eps,\xi}Z_i
		=\bigl({p_s^*} U_i^{{p_s^*}-1}-F_{\eps,\xi}\bigr)Z_i.
	\]
	For each fixed \(R>0\), uniformly on \(|t-\xi_i|\le R\),
	\[
		\begin{aligned}
			&V(t)=U_i(t)+O(\eps),\quad
			&&e^{\eps t}=1+O(\eps|\log\eps|),\\
			&V^{{p_s^*}+\eps-1}(t)=U_i^{{p_s^*}-1}(t)+O(\eps),\quad
			&&e^{-bt}=O(\eps).
		\end{aligned}
	\]
	Since \(U_i\) is bounded away from zero there, the power estimate also holds for \({p_s^*}-1<1\). Thus
	\begin{equation}\label{eq:LZi-local}
		\sup_{|t-\xi_i|\le R}|\mathfrak L_{\eps,\xi}Z_i(t)|
		\le C_R\eps|\log\eps|.
	\end{equation}
	All estimates are uniform for \(\Lambda\in\mathcal K_\delta\).
	Let \(c_Z\coloneqq\int_\R(U')^2>0\). Since
	\[
		\int_\R|Z_iZ_j|
		\le C\bigl(1+|\xi_i-\xi_j|\bigr)e^{-|\xi_i-\xi_j|}
		\quad(i\ne j),
	\]
	the Gram matrix satisfies
	\begin{equation}\label{eq:Gram}
		\int_\R Z_iZ_j=c_Z\delta_{ij}+O(\eps|\log\eps|).
	\end{equation}
	In particular, its inverse is uniformly bounded for small \(\eps\).

	\begin{lemma}[Bounds on the potential]\label{lem:F-bound}
		There are \(\eps_0>0\) and \(C>0\) such that, uniformly for
		\(0<\eps<\eps_0\) and \(\Lambda\in\mathcal K_\delta\),
		\[
			0<F_{\eps,\xi}(t)\le C\Wgt_{\theta_0}(t).
		\]
		More precisely,
		\begin{align*}
			({p_s^*}+\eps)e^{\eps t}V^{{p_s^*}+\eps-1}(t)
			&\le C\Wgt_{{p_s^*}-1}(t),\\
			p e^{-bt}V^{p-1}(t)
			&\le\begin{cases}
				C\eps e^{-(2p-{p_s^*}-1)(\xi_1-t)},&t\le\xi_1,\\
				C\eps\Wgt_{p-1}(t),&t\ge\xi_1.
			\end{cases}
		\end{align*}
		The same bounds hold for the first and second derivatives in \(\xi\). In particular,
		\[
			|D_\xi^jF_{\eps,\xi}(t)|\le C\Wgt_{\theta_0}(t),
			\qquad j=0,1,2.
		\]
	\end{lemma}

	\begin{proof}
		For \(t\le\xi_n\), \(e^{\eps t}V^\eps\le C\). For \(y\ge0\),
		\[
			e^{\eps(\xi_n+y)}V^\eps(\xi_n+y)\le Ce^{\eps\xi_n}\le C,
		\]
		and
		\[
			e^{-b(\xi_1-y)}V^{p-1}(\xi_1-y)
			\le Ce^{-b\xi_1}e^{-(p-1-b)y}
			\le C\eps e^{-(2p-{p_s^*}-1)y}.
		\]
		For \(t\ge\xi_1\), \(e^{-bt}\le C\eps\). The derivative bounds follow from the same estimates and \(|D_\xi V|+|D_\xi^2V|\le CV\).
	\end{proof}

\subsection{A preliminary weighted bound}

	The next lemma gives weighted decay for solutions in the energy space.

	\begin{lemma}[A priori weighted bound]\label{lem:apriori}
		Let \(h\) be measurable with \(\|h\|_\gamma<\infty\), and let
		\((\phi,c)\in H_k^s(\R)\times\R^n\) solve
		\eqref{eq:linear-projected}. Then \(\|\phi\|_\gamma<\infty\) and
		\begin{equation}\label{eq:green-rep-linear}
			\phi=G_s*g,
			\qquad g\coloneqq h+\sum_{i=1}^n c_iZ_i+F_{\eps,\xi}\phi.
		\end{equation}
	\end{lemma}

	\begin{proof}
		Since \(\|h\|_\gamma<\infty\), \(h\in L^2(\R)\cap L^\infty(\R)\). Lemma~\ref{lem:bounded-rep}, with \(\mathcal V=F_{\eps,\xi}\) and \(f_1=h+\sum_ic_iZ_i\), gives \eqref{eq:green-rep-linear} and a bounded continuous representative of \(\phi\). If \(|\phi|\le A\Wgt_\beta\) for some \(0<\beta\le\gamma\), then
		\[
			|h|\le\|h\|_\gamma\Wgt_\gamma,\qquad
			\sum_i|c_iZ_i|\le C|c|\Wgt_1,\qquad
			F_{\eps,\xi}|\phi|\le CA\Wgt_{\theta_0+\beta}.
		\]
		By Lemma~\ref{lem:comparison}, applied to these slowly varying majorants,
		\[
			|\phi(t)|\le C\int_\R e^{-|t-r|}
			\bigl(\|h\|_\gamma\Wgt_\gamma(r)+|c|\Wgt_1(r)+A\Wgt_{\theta_0+\beta}(r)\bigr)\,dr.
		\]
		At the initial step, \(|\phi|\le A_0\coloneqq\|\phi\|_{L^\infty}\), and the same argument uses \(A_0\Wgt_{\theta_0}\) in the last term.
		Convolution with \(e^{-|\cdot|}\) improves the decay exponent by
		\(\theta_0\) until it reaches \(\gamma\). The remaining terms are bounded by
		\(C\Wgt_\gamma\), since
		\[
			e^{-|\cdot|}*\Wgt_\gamma\le C\Wgt_\gamma,
			\qquad
			(1+|t-\xi_i|)e^{-|t-\xi_i|}\le C\Wgt_\gamma(t).
		\]
		For an intermediate exponent at least one, the convolution is bounded by \(C\sum_i(1+|t-\xi_i|)e^{-|t-\xi_i|}\le C\Wgt_\gamma\). Since \(\gamma<1\) and \(\theta_0>0\), finitely many steps suffice.
	\end{proof}

\subsection{Uniform invertibility}

	\begin{proposition}[Projected inverse]\label{prop:inverse}
		Assume \eqref{eq:p-range} and \eqref{eq:gamma-eta}, and fix \(n\) and
		\(\delta\). There are \(\eps_0>0\) and \(C>0\) with the following properties.
		\begin{enumerate}[label=\textup{(\alph*)}]
			\item For \(0<\eps<\eps_0\), \(\Lambda\in\mathcal K_\delta\) and
			measurable \(h\) with \(\|h\|_\gamma<\infty\), problem
			\eqref{eq:linear-projected} has a unique
			weak solution \((\phi,c)\in H_k^s(\R)\times\R^n\). It satisfies
			\[
				\|\phi\|_\gamma+\|\phi\|_{H_k^s}+\sum_{i=1}^n|c_i|
				\le C\|h\|_\gamma.
			\]
			\item If also \(h\in Y\), then \(\phi\in X\) and
			\[
				\|\phi\|_{\Lout}
				\le C\bigl(\|h\|_\gamma+\|h\|_{\eta,\Lout}\bigr)
				=C\|h\|_Y.
			\]
			\item For every \(\beta\in(0,\min\{1,2s\})\), \(\phi\in C^{0,\beta}(\R)\).
			If \(h\in C^{0,\sigma_0}(\R)\) for some \(\sigma_0>0\), then
			\(\phi\in C^{2s+\sigma}(\R)\) for some \(\sigma>0\), and the
			equation holds pointwise.
			\item Write \(T_{\eps,\xi}h\coloneqq\phi\) and let
			\(c_{\eps,\xi}\colon h\mapsto c\) denote the linear map giving the projection coefficients. Fix
			\(\xi^0=\xi(\Lambda^0)\), \(\Lambda^0\in\mathcal K_\delta\), and use the
			\(X\) and \(Y\) norms centered at \(\xi^0\). After decreasing \(\eps_0\), the map
			\[
				\xi\longmapsto\bigl(T_{\eps,\xi},c_{\eps,\xi}\bigr)
			\]
			is \(C^1\), in operator norm from \(Y\) to
			\((X\cap H_k^s)\times\R^n\), on the open set
			\(\mathcal N\coloneqq\{\xi:|\xi-\xi^0|<1\}\). Here
			\(X\cap H_k^s\) is equipped with the sum of the two norms. Moreover,
			\[
				\|\partial_{\xi_\ell}T_{\eps,\xi}h\|_X
				+\|\partial_{\xi_\ell}T_{\eps,\xi}h\|_{H_k^s}
				+|\partial_{\xi_\ell}c_{\eps,\xi}(h)|
				\le C\|h\|_Y.
			\]
			If \(\xi\mapsto h_\xi\) is \(C^1\) with values in \(Y\), then
			\[
				\partial_{\xi_\ell}(T_{\eps,\xi}h_\xi)
				=(\partial_{\xi_\ell}T_{\eps,\xi})h_\xi
				+T_{\eps,\xi}(\partial_{\xi_\ell}h_\xi).
			\]
		\end{enumerate}
	\end{proposition}

	\begin{proof}
		The proof follows \cite{dPDM2003,MR2381679}, using Lemmas~\ref{lem:indicial} and~\ref{lem:nondeg} in place of ODE arguments.

		\noindent\emph{Step 1.} The projection coefficients. Since \(Z_\ell\in H^s_k\) and \(\Ps Z_\ell\in L^2\), testing the equation with \(Z_\ell\) and using the symmetry of \(\Bform\) gives
		\[
			\sum_ic_i\int Z_iZ_\ell=\int\phi\,\mathfrak L_{\eps,\xi}Z_\ell-\int hZ_\ell .
		\]
		Since \(|Z_\ell|\le Ce^{-|t-\xi_\ell|}\), we have \(\left|\int hZ_\ell\right|\le C\|h\|_\gamma\). For \(R>1\), \eqref{eq:LZi-local} on \(|t-\xi_\ell|\le R\) and \(|\mathfrak L_{\eps,\xi}Z_\ell|\le Ce^{-|t-\xi_\ell|}\) elsewhere give
		\[
			\int\Wgt_\gamma|\mathfrak L_{\eps,\xi}Z_\ell|
			\le C_R\eps|\log\eps|+Ce^{-R}.
		\]
		Choose first \(R\) large and then \(\eps\) small, uniformly in \(\Lambda\in\mathcal K_\delta\). By \eqref{eq:Gram},
		\begin{equation}\label{eq:c-bound}
			|c|\le C\|h\|_\gamma+o(1)\|\phi\|_\gamma\qquad(\eps\downarrow0).
		\end{equation}

		\noindent\emph{Step 2.} The a priori estimate. By Lemma~\ref{lem:apriori}, every weak solution satisfies \(\|\phi\|_\gamma<\infty\). If \(\|\phi\|_\gamma\le C\|h\|_\gamma\) fails, there are \(\eps_j\to0\), \(\Lambda^j\in\mathcal K_\delta\) and solutions with \(\|\phi_j\|_\gamma=1\) 
		and \(\|h_j\|_\gamma\to0\). By \eqref{eq:c-bound}, \(c^j\to0\). Fix \(i\) and set \(\widehat\phi_j\coloneqq\phi_j(\cdot+\xi_i^j)\), so that 
		\(|\widehat\phi_j|\le n\) and
		\[
			g_j\coloneqq h_j+\sum_kc_k^jZ_k^j+F_j\phi_j
		\]
		is uniformly bounded. Proposition~\ref{prop:regularity}(i) gives 
		\[
			\|\widehat\phi_j\|_{C^{0,\beta}(\R)}\le C\|g_j\|_{L^\infty}
		\]
		for some \(\beta>0\), uniformly in \(j\). By the Arzel\`a--Ascoli theorem, a subsequence converges locally uniformly to \(\widehat\phi\) with \(|\widehat\phi|\le n\). On compact sets, \(V_j(\cdot+\xi_i^j)\to U\), 
		\(e^{\eps_j(\cdot+\xi_i^j)}\to1\) and \(e^{-b(\cdot+\xi_i^j)}=O(\eps_j)\), so \(F_j(\cdot+\xi_i^j)\to{p_s^*}U^{{p_s^*}-1}\) locally uniformly. For 
		\(\zeta\in C_c^\infty(\R)\) we have 
		\[
			\int\widehat\phi_j\Ps\zeta=\int g_j(\cdot+\xi_i^j)\zeta,
		\]
		and \(\Ps\zeta\in L^1\). Dominated convergence gives 
		\(\Ps\widehat\phi={p_s^*}U^{{p_s^*}-1}\widehat\phi\) in \(\mathcal D'(\R)\), and 
		\[
			\int\widehat\phi U'=\lim\int\phi_jZ_i^j=0.
		\] 
		By Lemma~\ref{lem:nondeg} and the orthogonality, \(\widehat\phi=cU'=0\). Thus \(\sup_{|t-\xi_i^j|\le R}|\phi_j(t)|\to0\) for every \(R\) and every \(i\). Let 
		\(\mathcal O_R\coloneqq\bigcup_i[\xi_i^j-R,\xi_i^j+R]\). On \(\mathcal O_R\) we have \(\Wgt_\gamma\ge e^{-\gamma R}\), and outside it \(F_j\le 
		Cne^{-\theta_0R}\) by Lemma~\ref{lem:F-bound}. Hence
		\[
			\|F_j\phi_j\|_\gamma\le Ce^{\gamma R}\sup_{\mathcal O_R}|\phi_j|+Ce^{-\theta_0R}.
		\]
		Let \(j\to\infty\) and then \(R\to\infty\). Thus \(\|F_j\phi_j\|_\gamma\to0\) and \(\|g_j\|_\gamma\to0\). Since \(\phi_j\in L^2\) and \(\Ps\phi_j=g_j\), Lemma~\ref{lem:resolvent} gives
		\[
			1=\|\phi_j\|_\gamma\le C\|g_j\|_\gamma\longrightarrow0,
		\]
		a contradiction. Testing with \(\phi\) and using the orthogonality,
		\[
			\|\phi\|_{H^s_k}^2=\int F\phi^2+\int h\phi\le C\bigl(\|\phi\|_\gamma^2+\|h\|_\gamma\|\phi\|_\gamma\bigr)\int\Wgt_\gamma^2\le C\|h\|_\gamma^2 .
		\]
		Together with \eqref{eq:c-bound}, this proves (a) and uniqueness.

		\noindent\emph{Step 3.} Existence. Let \(\mathscr H_\xi\coloneqq\{\phi\in H^s_k:\int\phi Z_i=0,\ i=1,\dots,n\}\). By \eqref{eq:Gram}, \(H^s_k=\mathscr 
		H_\xi\oplus\operatorname{span}\{Z_i\}\), and every continuous functional vanishing on \(\mathscr H_\xi\) is a unique combination of \(\psi\mapsto\int 
		Z_i\psi\). On \(\mathscr H_\xi\), with the inner product \(\Bform\), the projected equation reads
		\[
			\phi-\mathcal A\phi=\tilde h,
		\]
		where 
		\(\Bform(\mathcal A\phi,\psi)=\int F\phi\psi\) and \(\Bform(\tilde h,\psi)=\int h\psi\). If \(\phi_j\rightharpoonup0\), then \(\|F\phi_j\|_{L^2}\to0\), since \(F\to0\) at \(\pm\infty\) and \(H^s_k\hookrightarrow L^2_{\rm loc}\) is compact. Also, 
		\[
			\|\mathcal A\phi_j\|_{H^s_k}^2\le\|F\phi_j\|_{L^2}\|\mathcal A\phi_j\|_{L^2}\le C\|F\phi_j\|_{L^2}\|\mathcal A\phi_j\|_{H^s_k}.
		\]
		Thus \(\mathcal A\) is compact. Any element of \(\ker(I-\mathcal A)\) solves \eqref{eq:linear-projected} with \(h=0\), with coefficients determined by the residual functional, and vanishes by Step~2. The Fredholm alternative gives existence; the residual determines \(c\).

		\noindent\emph{Step 4.} The exterior estimate (b). Let \(M\coloneqq\|h\|_Y\). By Lemma~\ref{lem:apriori}, \(\phi=G_s*g\) with \(g=h+\sum c_iZ_i+F\phi\). For \(\lambda>0\), set
		\[
			m_\lambda(t)\coloneqq
			\begin{cases}
				e^{-\lambda(\xi_1-t)}&\text{if }t<\xi_1,\\
				n&\text{if }\xi_1\le t\le\xi_n,\\
				e^{-\lambda(t-\xi_n)}&\text{if }t>\xi_n.
			\end{cases}
		\]
		This function is \(ne^\lambda\)-slowly varying. By (a), \(|\phi|\le CMm_\gamma\). Suppose that \(|\phi|\le A_\lambda m_\lambda\) for some \(\gamma\le\lambda<1\). Then
		\[
			|h|\le Mm_\eta,\qquad
			\sum_i|c_iZ_i|\le CM\Wgt_1,\qquad
			F|\phi|\le CA_\lambda\Wgt_{\theta_0}m_\lambda.
		\]
		The majorants are slowly varying, uniformly in the centers and in \(\gamma\le\lambda<1\). By Lemma~\ref{lem:comparison},
		\[
			|\phi(t)|\le C\int_\R e^{-|t-r|}
			\bigl(Mm_\eta(r)+M\Wgt_1(r)+A_\lambda\Wgt_{\theta_0}(r)m_\lambda(r)\bigr)\,dr.
		\]
		Split the integrals at \(\xi_1\) and set \(t=\xi_1-y\), \(y\ge0\). By \eqref{eq:conv-elementary} and \(\eta>1\), the first two terms are bounded by \(CMe^{-y}\) and \(CM(1+y)e^{-y}\). For the last term, \(\Wgt_{\theta_0}m_\lambda\le n^2\) on \([\xi_1,\infty)\), while
		\[
			\Wgt_{\theta_0}(\xi_1-y')m_\lambda(\xi_1-y')\le n^2e^{-(\theta_0+\lambda)y'},\qquad y'\ge0.
		\]
		Thus
		\[
			|\phi(\xi_1-y)|\le C(M+A_\lambda)(1+y)e^{-y}
			+CA_\lambda\int_0^\infty e^{-|y-y'|}e^{-(\theta_0+\lambda)y'}\,dy'.
		\]
		The same estimate holds at \(\xi_n+y\). Choose \(\theta'\in(0,\theta_0]\) such that \(\gamma+k\theta'\ne1\) for every integer \(k\ge0\), and let \(k_0\) be the first integer with \(\gamma+k_0\theta'>1\). Since \(\theta'\le\theta_0\), \eqref{eq:conv-elementary} and \((1+y)e^{-y}\le C_\beta e^{-\beta y}\) for \(\beta<1\) show that the exterior exponent improves from \(\gamma+k\theta'\) to \(\gamma+(k+1)\theta'\) whenever the latter is below one. Iterate on both half-lines, retaining the interior bound from (a); at each step \(A_\lambda=C_\lambda M\). In the final step,
		\[
			\lambda=\gamma+(k_0-1)\theta'<1,
			\qquad \theta_0+\lambda\ge\gamma+k_0\theta'>1,
		\]
		so the integral term is bounded by \(CMe^{-y}\). Thus
		\[
			|\phi(\xi_1-y)|+|\phi(\xi_n+y)|\le CM(1+y)e^{-y}.
		\]
		For a source \(h+h_1\) with \(\|h_1\|_1<\infty\), replace \(M\) by \(\|h\|_Y+\|h_1\|_1\). Part (a) applies since \(\|h_1\|_\gamma\le\|h_1\|_1\). The added term satisfies \(|h_1|\le\|h_1\|_1\Wgt_1\) and
		\[
			|(G_s*h_1)(t)|\le C\|h_1\|_1\sum_i(1+|t-\xi_i|)e^{-|t-\xi_i|}.
		\]

		\noindent\emph{Step 5.} Regularity. By Lemma~\ref{lem:apriori}, \(\phi\) is bounded and
		\[
			\Ps\phi=g\coloneqq h+\sum_ic_iZ_i+F\phi
		\]
		in \(\mathcal D'(\R)\), with \(g\in L^\infty(\R)\). Proposition~\ref{prop:regularity}(i) gives the first assertion. The functions \(F\) and \(Z_i\) are bounded and Lipschitz: Lemma~\ref{lem:F-bound} bounds \(F\), and \(|V'|\le CV\), together with the factors \(e^{\eps t}\) and \(e^{-bt}\), gives \(|\partial_tF|\le CF\). Thus, if \(h\in C^{0,\sigma_0}(\R)\), then \(g\in C^{0,\min\{\sigma_0,\beta\}}(\R)\). Proposition~\ref{prop:regularity}(ii), with \(0<\sigma<\min\{\sigma_0,\beta,1,2-2s\}\) and \(2s+\sigma\ne1\), gives the second assertion.

		\noindent\emph{Step 6.} Dependence on the centers. By \eqref{eq:xigap}, \(|\log\Lambda_i-\log\Lambda_i^0|<2\) on \(\mathcal N\). Replace \(\delta\) by \(\delta e^{-2}\) and use the norms centered at \(\xi^0\) throughout \(\mathcal N\), as allowed by Lemma~\ref{lem:moving-norms}.
		
		\begin{enumerate}[label=\textup{(\roman*)}]
			\item  For \(h\in Y\), \(\|h_1\|_1<\infty\) and \(\mathfrak c\in\R^n\), the system \(\mathfrak L_{\eps,\xi}w=h+h_1+\sum_id_iZ_i\), \(\int wZ_i=\mathfrak c_i\), has a unique solution, and
			\[
				\|w\|_X+\|w\|_{H^s_k}+|d|\le C\bigl(\|h\|_Y+\|h_1\|_1+|\mathfrak c|\bigr).
			\]
			Write \(w=\tilde w+\sum_ja_jZ_j\), with \(a\) given by the Gram system for \(\mathfrak c\). Then \(\tilde w\) solves \eqref{eq:linear-projected} with source \(h+h_1-\sum_ja_j\mathfrak L_{\eps,\xi}Z_j\). Since \(|\mathfrak L_{\eps,\xi}Z_j(t_\pm(y))|\le Ce^{-(1+\theta_0)y}\) and \(\eta<1+\theta_0\), \(\|\mathfrak L_{\eps,\xi}Z_j\|_Y\le C\). Steps~1--4 give the estimate.
			
			\item By Lemma~\ref{lem:F-bound} and Taylor's formula,
			\[
				\|(F_{\xi'}-F_\xi)w\|_Y\le C|\xi'-\xi|\|w\|_X,
			\]
			with a remainder \(O(|\xi'-\xi|^2)\|w\|_X\) after subtracting \(D_\xi F_\xi[\xi'-\xi]w\). Moreover,
			\[
				Z_i^{\xi'}-Z_i^\xi=-(\xi_i'-\xi_i)U_i''+O(|\xi'-\xi|^2e^{-|t-\xi_i|})
			\]
			in the norm \(\|\cdot\|_1\).
			
			\item Let \((\phi,c)\) and \((\phi',c')\) correspond to \(\xi\) and \(\xi'\). The difference \(\phi'-\phi\) solves the system in (i) at \(\xi\), with source \((F_{\xi'}-F_\xi)\phi'+\sum_ic_i'(Z_i^{\xi'}-Z_i^\xi)\) and constraints \(-\int\phi'(Z_i^{\xi'}-Z_i^\xi)\). By (i)--(ii),
			\[
				\|\phi'-\phi\|_X+\|\phi'-\phi\|_{H^s_k}+|c'-c|\le C|\xi'-\xi|\|h\|_Y.
			\]
			
			\item Let \((\dot\phi,\dot c)\) solve the system in (i) with source \((\partial_{\xi_\ell}F)\phi+c_\ell\,\partial_{\xi_\ell}Z_\ell\) and constraints \(-\delta_{i\ell}\int\phi\,\partial_{\xi_\ell}Z_\ell\). For \(\xi'=\xi+\tau e_\ell\), (ii)--(iii) give data of size \(O(|\tau|\|h\|_Y)\) for the equation satisfied by the difference quotient minus \((\dot\phi,\dot c)\). Part (i) gives differentiability and the bound in (d). The derivative is continuous in operator norm: \(|D^2_\xi F|\le C\Wgt_{\theta_0}\) by Lemma~\ref{lem:F-bound}, and \(\eta<1+\theta_0\), so \(\partial_{\xi_\ell}F_\xi\) is Lipschitz as a multiplier from \(X\) to \(Y\); the argument in (iii), with the constraint data included, gives Lipschitz dependence of the solution operator in (i); and \(c_{\eps,\xi}\) and \(\partial_{\xi_\ell}Z_\ell\) are Lipschitz, the latter in \(\|\cdot\|_1\). The last formula follows from the chain rule.
			
		\end{enumerate}
	\end{proof}

\section{The projected nonlinear problem}\label{sec:nonlinear}

	Recall that \(q={p_s^*}+\eps\). Set
	\begin{align*}
		\mathfrak R_{\eps,\xi}
		&\coloneqq e^{\eps t}V^{{p_s^*}+\eps}+e^{-bt}V^p-\sum_{i=1}^nU_i^{p_s^*},
		\\
		\mathfrak N_{\eps,\xi}(\phi)
		&\coloneqq e^{\eps t}\Big[(V+\phi)_+^{{p_s^*}+\eps}-V^{{p_s^*}+\eps}
		-({p_s^*}+\eps)V^{{p_s^*}+\eps-1}\phi\Big]\\
		&\quad +e^{-bt}\Big[(V+\phi)_+^p-V^p-pV^{p-1}\phi\Big].
	\end{align*}
	Since \(\Ps V=\sum_iU_i^{p_s^*}\),
	\[
		F_{\eps,\xi}\phi+\mathfrak R_{\eps,\xi}
		+\mathfrak N_{\eps,\xi}(\phi)
		=e^{\eps t}(V+\phi)_+^q+e^{-bt}(V+\phi)_+^p-\sum_iU_i^{p_s^*}.
	\]
	We seek \((\phi,c)\) such that
	\begin{equation}\label{eq:projected-nl}
		\begin{cases}
			\mathfrak L_{\eps,\xi}\phi
			=\mathfrak N_{\eps,\xi}(\phi)+\mathfrak R_{\eps,\xi}
			+\displaystyle\sum_{i=1}^n c_iZ_i,\\[2mm]
			\displaystyle\int_{\R}\phi Z_i=0,
			\qquad i=1,\ldots,n.
		\end{cases}
	\end{equation}
	Constants may depend on \(N,s,p,n,\delta,\gamma,\eta\) and \(\alpha\), and are uniform for \(\Lambda\in\mathcal K_\delta\) and small \(\eps>0\). Parameter derivatives are taken at fixed \(t\), in the fixed spaces of Lemma~\ref{lem:moving-norms} and Proposition~\ref{prop:inverse}(d).

\subsection{Remainder estimates}

	\begin{lemma}[Estimate for the approximation error]\label{lem:R}
		Choose \(\gamma,\eta\) as in \eqref{eq:gamma-eta} and let \(\alpha\) be as in \eqref{eq:gamma-alpha-thm}. Then
		\[
			\|\mathfrak R_{\eps,\xi}\|_\gamma
			+\|D_\xi\mathfrak R_{\eps,\xi}\|_\gamma
			\le C\eps^\alpha,\qquad
			\|\mathfrak R_{\eps,\xi}\|_{\eta,\Lout}
			+\|D_\xi\mathfrak R_{\eps,\xi}\|_{\eta,\Lout}
			\le C\eps^\alpha.
		\]
		The map \(\xi\mapsto\mathfrak R_{\eps,\xi}\) is \(C^1\) with values in \(Y\).
	\end{lemma}

	\begin{proof}
		Write \(\mathfrak R=R_1+R_2+R_3+R_4\), where
		\[
			R_1\coloneqq e^{\eps t}(V^{{p_s^*}+\eps}-V^{p_s^*}),\ \
			R_2\coloneqq(e^{\eps t}-1)V^{p_s^*},\ \
			R_3\coloneqq V^{p_s^*}-\textstyle\sum_iU_i^{p_s^*},\ \
			R_4\coloneqq e^{-bt}V^p.
		\]
		Fix a small \(\iota>0\) with \(2\iota<{p_s^*}-\gamma\) and \(\iota<{p_s^*}-\eta\). The bounds 
		\[
			|V^{{p_s^*}+\eps}-V^{p_s^*}|\le C\eps V^{{p_s^*}-\iota}\quad\text{ and }\quad 
			|e^{\eps t}-1|\le C\eps\bigl(|\log\eps|+\dist(t,\{\xi_i\})\bigr)e^{\eps(t-\xi_n)_+}
		\] 
		and \(\eps\xi_n=o(1)\), together with \(V^{p_s^*}\le Ce^{-{p_s^*}y}\) on the exterior half-lines, give \(\|R_1\|_Y+\|R_2\|_Y\le C\eps|\log\eps|\). On the cell of \(\xi_i\), \(|R_3|\le CU_i^{{p_s^*}-1}\sum_{j\ne i}U_j\). For centers at distance \(L\),
		\[
			\sup_{t\in\Omega_i} U_i^{{p_s^*}-1}U_j/\Wgt_\gamma\le C\max\{e^{-L},e^{-({p_s^*}-\gamma)L/2}\},
		\]
		and \(e^{-L}\le C\eps\) gives \(O(\eps^{\min\{1,({p_s^*}-\gamma)/2\}})\). On the left half-line, \(U_1\le Ce^{-y}\) and \(\sum_{j\ge2}U_j\le C\eps e^{-y}\), hence 
		\[
			|R_3|\le CU_1^{{p_s^*}-1}\sum_{j\ge2}U_j\le C\eps e^{-{p_s^*}y}.
		\] 
		The right half-line is symmetric. On the left half-line 
		\[
			R_4\le Ce^{-b\xi_1}e^{by}e^{-py}=C\eps e^{-(2p-{p_s^*})y},
		\] 
		with \(2p-{p_s^*}>\eta\). On the right half-line, 
		\[
			R_4\le Ce^{-b\xi_n}e^{-(p+b)y}\le C\eps e^{-{p_s^*}y},\text{ and  }e^{-bt}\le C\eps\text{ on }[\xi_1,\xi_n].
		\] 
		Differentiation in \(\xi\) gives 
		\begin{align*}
			&\partial_{\xi_\ell}R_1=-e^{\eps t}\bigl[({p_s^*}+\eps)V^{{p_s^*}+\eps-1}-{p_s^*}V^{{p_s^*}-1}\bigr]Z_\ell,
			&&\partial_{\xi_\ell}R_2=-{p_s^*}(e^{\eps t}-1)V^{{p_s^*}-1}Z_\ell,\\
			&\partial_{\xi_\ell}R_3=-{p_s^*}(V^{{p_s^*}-1}-U_\ell^{{p_s^*}-1})Z_\ell,
			&&\partial_{\xi_\ell}R_4=-pe^{-bt}V^{p-1}Z_\ell.
		\end{align*}
		Since \(|Z_\ell|\le CV\), the derivatives of \(R_1,R_2,R_4\) satisfy the same bounds. For \(R_3\), the mean value theorem on \([U_\ell,CU_\ell]\) in the cell of \(\xi_\ell\), and \(V\le CU_i\) in the other cells, give
		\[
			|V^{{p_s^*}-1}-U_\ell^{{p_s^*}-1}||Z_\ell|
			\le
			\begin{cases}
				CU_\ell^{{p_s^*}-1}\displaystyle\sum_{j\ne\ell}U_j,&\text{in the cell of }\xi_\ell,\\
				CU_i^{{p_s^*}-1}U_\ell,&\text{in the cell of }\xi_i,\ i\ne\ell.
			\end{cases}
		\]
		Thus \(\partial_{\xi_\ell}R_3\) satisfies the bounds of \(R_3\), including when \({p_s^*}<2\), and
		\[
			|\partial_{\xi_\ell}R_3(t_\pm(y))|\le C\eps e^{-{p_s^*}y}.
		\]
		For second derivatives, use \(|\partial_\xi V|+|\partial^2_\xi V|\le CV\) and the mean value theorem on \([U_i,CU_i]\) in the cell of \(\xi_i\):
		\[
			|\partial^2_\xi R_3|\le CU_i^{{p_s^*}-1}\sum_{j\ne i}U_j.
		\]
		These bounds also hold when \({p_s^*}<2\) and give the \(C^1\) dependence.
	\end{proof}

\subsection{Nonlinear estimates and the fixed point argument}

	For \(r>0\), \(A>0\) and \(B\in\R\), set
	\[
		Q_r(A,B)\coloneqq(A+B)_+^r-A^r-rA^{r-1}B.
	\]
	The elementary bounds are
	\begin{equation}\label{eq:Qr-estimates}
		|Q_r(A,B)|\le
		\begin{cases}
			C\bigl(|B|^r+A^{r-1}|B|\bigr),&0<r<1,\\
			C|B|^r,&1\le r\le2,\\
			C\bigl(A^{r-2}B^2+|B|^r\bigr),&r\ge2.
		\end{cases}
	\end{equation}
	For \(r>1\) and \(M\coloneqq\max\{|B_1|,|B_2|\}\),
	\begin{equation}\label{eq:Qr-lipschitz}
		\begin{aligned}
			|Q_r(A,B_1)-Q_r(A,B_2)|
			&\le C\Xi_r|B_1-B_2|,\\
			\Xi_r&=
			\begin{cases}
				M^{r-1},&1<r\le2,\\
				A^{r-2}M+M^{r-1},&r\ge2.
			\end{cases}
		\end{aligned}
	\end{equation}
	For \(1<r\le2\), these follow from the \((r-1)\)-H\"older continuity of \(x\mapsto x_+^{r-1}\). For \(r\ge2\), integrate the derivative. For \(0<r<1\), use \(|(A+B)_+^r-A^r|\le|B|^r\). For \(r=1\),
	\[
		Q_1(A,B)=(-A-B)_+,
		\qquad |Q_1(A,B)|\le |B|,
	\]
	because \(A>0\).

	\begin{lemma}[Nonlinear estimates]\label{lem:N}
		Let \(\vartheta\coloneqq\min\{1,p-1\}>0\). There are \(\mathfrak r_0>0\) and \(C>0\) such that, if
		\begin{equation}\label{eq:small-two-norms}
			\|\phi\|_X\le\mathfrak r\le\mathfrak r_0,
		\end{equation}
		then
		\[
			\|\mathfrak N_{\eps,\xi}(\phi)\|_\gamma
			\le C\mathfrak r^\vartheta\|\phi\|_X,
			\qquad
			\|\mathfrak N_{\eps,\xi}(\phi)\|_{\eta,\Lout}
			\le C\mathfrak r^\vartheta\|\phi\|_X.
		\]
		If \(\phi_1,\phi_2\) satisfy \eqref{eq:small-two-norms}, then
		\[
			\|\mathfrak N(\phi_1)-\mathfrak N(\phi_2)\|_Y
			\le C\mathfrak r^\vartheta\|\phi_1-\phi_2\|_X.
		\]
		The map \((\xi,\phi)\mapsto\mathfrak N_{\eps,\xi}(\phi)\) is \(C^1\) on the corresponding open subset of the fixed space \(\R^n\times X\), with values in \(Y\). Moreover,
		\[
			\begin{aligned}
				\partial_{\xi_\ell}\mathfrak N_{\eps,\xi}(\phi)
				&=-q e^{\eps t}Q_{q-1}(V,\phi)Z_\ell
				-p e^{-bt}Q_{p-1}(V,\phi)Z_\ell,\\
				\|\partial_{\xi_\ell}\mathfrak N_{\eps,\xi}(\phi)\|_Y
				&\le C\mathfrak r^\vartheta.
			\end{aligned}
		\]
	\end{lemma}

	\begin{proof}
		By definition,
		\[
			\mathfrak N(\phi)=e^{\eps t}Q_q(V,\phi)+e^{-bt}Q_p(V,\phi).
		\]
		Apply \eqref{eq:Qr-estimates}. On \([\xi_1,\xi_n]\), \(|\phi|^r\le C\|\phi\|_\gamma^r\Wgt_\gamma\) for \(r>1\), and \(e^{-bt}\le C\eps\). On the exterior half-lines, \(|\phi|\le\|\phi\|_X(1+y)e^{-y}\). The weighted terms decay with exponents \(2p-{p_s^*}\) and \({p_s^*}\), both greater than \(\eta\), so the polynomial factors are absorbed. The Lipschitz bound follows from \eqref{eq:Qr-lipschitz}. The derivative formulas follow from \(\partial_AQ_r=rQ_{r-1}\); continuity follows from the H\"older continuity of \(x_+^{r-1}\).
	\end{proof}

	\begin{proposition}[Solution of the projected nonlinear problem]\label{prop:nonlinear}
		Let \(\alpha\) be as in \eqref{eq:gamma-alpha-thm}. There are \(\eps_0>0\) and \(C_0>0\) such that, for \(\Lambda\in\mathcal K_\delta\) and \(0<\eps<\eps_0\), problem \eqref{eq:projected-nl} has a unique solution \((\phi,c)\in(X\cap H^s_k(\R))\times\R^n\) with
		\[
			\phi\in\mathcal B_\eps
			\coloneqq\{\psi\in X:\|\psi\|_X\le C_0\eps^\alpha\}.
		\]
		Moreover,
		\[
			\|\phi\|_{H^s_k}+\sum_i|c_i|\le C\eps^\alpha.
		\]
		For every \(0<\sigma<\min\{1,2-2s\}\) with \(2s+\sigma\ne1\), \(\phi\in C^{2s+\sigma}(\R)\), and the equation holds pointwise. These conclusions also hold, with the same \(\eps_0\) and \(C_0\), for \(\Lambda\in\mathcal K_{\delta e^{-2}}\). In the fixed spaces of Proposition~\ref{prop:inverse}(d), the map \(\xi\mapsto(\phi(\xi),c(\xi))\) is defined and \(C^1\) into \((X\cap H^s_k(\R))\times\R^n\) on the open set \(\{\xi:|\xi-\xi(\Lambda)|<1\text{ for some }\Lambda\in\mathcal K_\delta\}\).
	\end{proposition}

	\begin{proof}
		Proposition~\ref{prop:inverse} and Lemmas~\ref{lem:R}--\ref{lem:N} hold for every \(\delta\in(0,1)\); apply them with \(\delta e^{-2}\) in place of \(\delta\), and choose \(\eps_0\) and \(C_0\) uniformly on \(\mathcal K_{\delta e^{-2}}\). The map
		\[
			\phi\mapsto T_{\eps,\xi}(\mathfrak R+\mathfrak N(\phi))
		\]
		maps \(\mathcal B_\eps\) into itself and is a contraction for \(C_0\) large and \(\eps\) small. At the fixed point,
		\[
			\|\phi\|_X\le C\|\mathfrak R\|_Y+C(C_0\eps^\alpha)^\vartheta\|\phi\|_X
			\le \frac{C_0}{2}\eps^\alpha,
		\]
		after increasing \(C_0\) and decreasing \(\eps_0\). Thus \(\phi\in\operatorname{int}\mathcal B_\eps\). For any fixed \(K\ge1\), the same argument gives contraction and uniqueness in the ball of radius \(KC_0\eps^\alpha\), after decreasing \(\eps_0\). By Lemma~\ref{lem:moving-norms}, the balls for two configurations at distance at most one are comparable. Uniqueness therefore makes \(\phi\) independent of the weights used in the construction.

		For regularity, set \(v\coloneqq V+\phi\). Then
		\[
			\Ps\phi=g\coloneqq e^{\eps t}(v_+)^q+e^{-bt}(v_+)^p-\sum_iU_i^{p_s^*}+\sum_ic_iZ_i
		\]
		in \(\mathcal D'(\R)\). On the exterior half-lines, \(\phi\in X\) gives \(|v|\le C(1+y)e^{-y}\), where \(y=\dist(t,[\xi_1,\xi_n])\). Thus \(e^{-bt}|v|^{p-1}\) and \(e^{\eps t}|v|^{q-1}\) are bounded, and \(g\in L^\infty(\R)\). Proposition~\ref{prop:regularity}(i) gives \(\phi\in C^{0,\beta}(\R)\) for every \(\beta\in(0,\min\{1,2s\})\). Use the majorant
		\[
			|v(t)|\le\mathcal M(t)\coloneqq
			\begin{cases}
				C\Wgt_1(t)+C_0\eps^\alpha\Wgt_\gamma(t),&t\in[\xi_1,\xi_n],\\
				C(1+y)e^{-y},&t\notin[\xi_1,\xi_n].
			\end{cases}
		\]
		It is bounded and slowly varying, and \(e^{-bt}\mathcal M^{p-1}\), \(e^{\eps t}\mathcal M^{q-1}\) are bounded. Since \(V,U_i^{p_s^*},Z_i\) are bounded and Lipschitz, Lemma~\ref{lem:holder-transfer}, with \(r\in\{p,q\}\), transfers the global H\"older regularity of \(\phi\) to \(g\), up to order one. Finitely many applications of Proposition~\ref{prop:regularity}(ii) give the assertion.

		For the \(C^1\) dependence, fix \(\bar\xi\) with \(\Lambda(\bar\xi)\in\mathcal K_\delta\) and take \(\xi^0=\bar\xi\). By \eqref{eq:xigap}, \(|\xi-\bar\xi|<1\) implies \(\Lambda(\xi)\in\operatorname{int}\mathcal K_{\delta e^{-2}}\), so the construction applies throughout \(\mathcal N\). In the fixed spaces of Proposition~\ref{prop:inverse}(d), decrease \(\mathfrak r_0\) using Lemma~\ref{lem:moving-norms} and set \(\mathcal G(\xi,\psi)\coloneqq\psi-T_{\eps,\xi}(\mathfrak R+\mathfrak N(\psi))\) on \(\mathcal N\times\{\psi\in X:\|\psi\|_X<\mathfrak r_0\}\). Proposition~\ref{prop:inverse}(d) and Lemma~\ref{lem:N} give
		\[
			\|D_\psi\mathcal G(\xi,\psi)-I\|_{\mathcal L(X)}
			\le C\|\psi\|_X^\vartheta.
		\]
		At each fixed point, \(\|\phi(\xi)\|_X=O(\eps^\alpha)\), so \(D_\psi\mathcal G(\xi,\phi(\xi))=I+O(\eps^{\alpha\vartheta})\) is invertible for small \(\eps\). Fix \(\xi_*\in\mathcal N\). The implicit function theorem \cite[Chapter~2]{AmbrosettiProdi1993} gives a \(C^1\) branch \(\tilde\phi\) near \(\xi_*\), with \(\tilde\phi(\xi_*)=\phi(\xi_*)\). By continuity and Lemma~\ref{lem:moving-norms},
		\[
			\|\tilde\phi(\xi)\|_{X,\xi}
			\le C\|\tilde\phi(\xi)\|_{X,\bar\xi}
			\le KC_0\eps^\alpha
		\]
		near \(\xi_*\), with \(K\) depending only on \(n\), \(\gamma\) and \(\eta\). Uniqueness in this ball gives \(\tilde\phi=\phi\). Since \(\xi_*\) is arbitrary, \(\phi\) is \(C^1\) throughout \(\mathcal N\). Finally,
		\[
			(\phi,c)=(T_{\eps,\xi},c_{\eps,\xi})(\mathfrak R+\mathfrak N(\phi))
		\]
		and Proposition~\ref{prop:inverse}(d) give \(C^1\) dependence into \((X\cap H^s_k(\R))\times\R^n\).
	\end{proof}

\subsection{Positivity and behavior at infinity}

	When \(c=0\), the exterior estimate improves to pure exponential decay.

	\begin{lemma}[Positivity and regularity in the original variables]\label{lem:positive}
		Let \((\phi,c)\) be the solution of Proposition~\ref{prop:nonlinear} with \(\Lambda\in\mathcal K_\delta\), and assume \(c_1=\cdots=c_n=0\). Then, for \(\eps\) small, the following hold.
		\begin{enumerate}[label=\textup{(\alph*)}]
			\item For \(y\ge0\),
			\begin{equation}\label{eq:phi-tail-improved}
				|\phi(\xi_1-y)|+|\phi(\xi_n+y)|
				\le C\eps^\alpha e^{-y}.
			\end{equation}
			\item With \(\varrho\) as in Corollary~\ref{cor:asymptotics},
			\begin{equation}\label{eq:relative-error-EF}
				\sup_{t\in\R}\frac{|\phi(t)|}{V(t)}\le C\eps^\varrho.
			\end{equation}
			For \(n=1\), the right-hand side can be replaced by \(C\eps^\alpha\). In particular, \(v=V+\phi>0\) and
			\begin{equation}\label{eq:v-tail}
				c\Wgt_1(t)\le v(t)\le C\Wgt_1(t),\qquad t\in\R.
			\end{equation}
			\item The function
			\[
				u(r)=r^{-a}v(-a\log r),\qquad r>0,
			\]
			extends to a bounded positive radial solution of \eqref{eq:main-problem} on \(\R^N\). It belongs to \(C^{2s+\sigma}(\R^N)\) for every \(0<\sigma<\min\{1,2-2s\}\) with \(2s+\sigma\ne1\), and satisfies
			\begin{equation}\label{eq:u-decay}
				C^{-1}\mu_1^a|x|^{2s-N}\le u(x)\le C\mu_1^a|x|^{2s-N},
				\qquad |x|\ge\mu_1,
				\qquad \mu_1=e^{-\xi_1/a}.
			\end{equation}
		\end{enumerate}
	\end{lemma}

	\begin{proof}
		(a) Since \(c=0\), \(\phi=G_s*g\) with \(g=F\phi+\mathfrak R+\mathfrak N(\phi)\). The exterior bounds and \(\eta<1+\theta_0\) give
		\begin{align*}
			F(t_\pm(y))(1+y)e^{-y}
			&\le C(1+y)e^{-(1+\theta_0)y}\le Ce^{-\eta y},\\
			|g(t_\pm(y))|&\le C\eps^\alpha e^{-\eta y}.
		\end{align*}
		On \([\xi_1,\xi_n]\), \(|g|\le Cn\eps^\alpha\), whose contribution at \(t_\pm(y)\) is at most \(C\eps^\alpha e^{-y}\). Since \(\eta>1\), Lemma~\ref{lem:comparison} and \eqref{eq:conv-elementary} give \eqref{eq:phi-tail-improved}.

		(b) On \([\xi_1,\xi_n]\),
		\begin{align*}
			\frac{|\phi|}{V}&\le C\eps^\alpha\frac{\Wgt_\gamma}{\Wgt_1}
			\le C\eps^\alpha e^{(1-\gamma)D(t)},\\
			D(t)&=\dist(t,\{\xi_i\})\le\frac{|\log\eps|+C}{2}.
		\end{align*}
		On the exterior half-lines, (a) gives \(|\phi|/V\le C\eps^\alpha\).

		(c) The function \(u(r)=e^tv(t)\) is comparable to \(e^{\xi_n}\) near \(0\) and to \(\mu_1^ar^{2s-N}\) for \(r\ge\mu_1\). By Proposition~\ref{prop:EF} and Lemma~\ref{lem:removable}, it solves the equation in \(\R^N\setminus\{0\}\) and across \(0\), respectively. Also, \(v\in L^2(\R)\) and
		\[
			\Ps v=H\coloneqq e^{-bt}v^p+e^{\eps t}v^{{p_s^*}+\eps}\in L^2(\R),
		\]
		so \(v=G_s*H\) by Lemma~\ref{lem:indicial}(iv). By \eqref{eq:Riesz-cylinder},
		\[
			u=c^{\mathrm R}_{N,s}|\cdot|^{2s-N}*f,
		\]
		with \(f\coloneqq u^p+u^{{p_s^*}+\eps}\). Boundedness of \(u\), the decay \(u\le C|x|^{2s-N}\) for \(|x|\ge\mu_1\), and \(p(N-2s)>N\) give \(f\in L^1(\R^N)\cap L^\infty(\R^N)\). Choose \(\beta_0\in(0,\min\{1,2s\})\) with \(\beta_0+2sj\ne1\) for every integer \(j\ge0\). By Lemma~\ref{lem:Riesz}, \(u\in C^{0,\beta_0}(\R^N)\). Since \(p,q>1\), \(f\) inherits the H\"older regularity of \(u\), up to order one. Iterating \eqref{eq:Silvestre-holder}, with exponents below \(\min\{1,2-2s\}\), gives Lipschitz regularity of \(u\) and \(f\) after finitely many steps. A final application gives \(u\in C^{2s+\sigma}(\R^N)\) for every \(0<\sigma<\min\{1,2-2s\}\) with \(2s+\sigma\ne1\).
	\end{proof}

\section{The reduced energy}\label{sec:energy}

	We expand the energy in terms of the centers, following the local tower constructions \cite{dPDM2003,MR2381679}. The cylindrical kernel determines the interactions. With the truncation of the projected equation, set
	\[
		\E_\eps(v)
		\coloneqq\frac12\Bform(v,v)
		-\frac{1}{{p_s^*}+\eps+1}\int_{\R}e^{\eps t}(v_+)^{{p_s^*}+\eps+1}\,dt
		-\frac{1}{p+1}\int_{\R}e^{-bt}(v_+)^{p+1}\,dt .
	\]
	Positivity of \(V+\phi\) is proved only when the projection coefficients vanish, so the positive parts are retained. Set
	\[
		\mathcal Y\coloneqq H^s_k(\R)\cap X,
		\qquad
		\|v\|_{\mathcal Y}\coloneqq\|v\|_{H^s_k}+\|v\|_X .
	\]
	Fix the weights at a reference configuration, as in Lemma~\ref{lem:moving-norms}.

	For \(v\) in bounded subsets of \(X\), \(|v|\le C\mathcal M_\xi\), where
	\[
		\mathcal M_\xi(t)=
		\begin{cases}
			(1+\xi_1-t)e^{t-\xi_1},&t<\xi_1,\\
			\Wgt_\gamma(t),&\xi_1\le t\le\xi_n,\\
			(1+t-\xi_n)e^{\xi_n-t},&t>\xi_n.
		\end{cases}
	\]
	Since \(p+1-b>0\),
	\[
		\int_\R e^{-bt}\mathcal M_\xi^{p+1}<\infty,\qquad
		\int_\R e^{\eps t}\mathcal M_\xi^{q+1}<\infty,\qquad q={p_s^*}+\eps.
	\]
	Dominated convergence gives \(C^2\) regularity of the nonlinear terms on \(X\). With \(\vartheta=\min\{1,p-1\}\) as in Lemma~\ref{lem:N}, the local \(\vartheta\)-H\"older continuity of \(x_+^{p-1}\) and \(x_+^{q-1}\), since \(1<p<{p_s^*}<q\), gives, for \(v,w\in\mathcal Y\) with bounded \(X\)-norms and \(\psi,\chi\in\mathcal Y\), with \(D^2\E_\eps\) given by \eqref{eq:energy-derivatives} below,
	\[
		\bigl|(D^2\E_\eps(v)-D^2\E_\eps(w))[\psi,\chi]\bigr|
		\le C\|v-w\|_X^{\vartheta}\|\psi\|_X\|\chi\|_X.
	\]
	Since \(\Bform\) cancels in this difference, the left-hand side involves only the nonlinear terms, and the same bound holds for \(v,w\) in bounded subsets of \(X\) and \(\psi,\chi\in X\). Hence these terms are \(C^{2,\vartheta}\) on bounded subsets of \(X\); see \cite[Ch.~1, \S2]{AmbrosettiProdi1993} for Nemytskii operators on Lebesgue spaces.

	The continuity of \(\Bform\) on \(H^s_k\) gives \(\E_\eps\in C^2(\mathcal Y)\), with
	\begin{equation}\label{eq:energy-derivatives}
		\begin{aligned}
			D\E_\eps(v)[\psi]
			&=\Bform(v,\psi)
			-\int_\R\bigl(e^{\eps t}(v_+)^{{p_s^*}+\eps}
			+e^{-bt}(v_+)^p\bigr)\psi\,dt,\\
			D^2\E_\eps(v)[\psi,\chi]
			&=\Bform(\psi,\chi)
			-\int_\R\bigl(({p_s^*}+\eps)e^{\eps t}(v_+)^{{p_s^*}+\eps-1}
			+p e^{-bt}(v_+)^{p-1}\bigr)\psi\chi\,dt.
		\end{aligned}
	\end{equation}
	By Proposition~\ref{prop:H-equivalence}(iii), \(V\in\mathcal Y\); Proposition~\ref{prop:nonlinear} gives \(\phi\in\mathcal Y\), hence \(V+\phi\in\mathcal Y\). Since \(V>0\), the truncation does not
	affect \(\E_\eps(V)\).

	Using \(\Bform(U,U)=\int_\R U^{{p_s^*}+1}\), define
	\begin{equation}\label{eq:a3}
		\begin{aligned}
			a_1&\coloneqq\left(\frac12-\frac1{{p_s^*}+1}\right)
			\int_\R U^{{p_s^*}+1}>0,\qquad
			a_3\coloneqq\frac1{{p_s^*}+1}\int_\R U^{{p_s^*}+1}>0,\\
			a_4&\coloneqq\frac1{({p_s^*}+1)^2}\int_\R U^{{p_s^*}+1}
			-\frac1{{p_s^*}+1}\int_\R U^{{p_s^*}+1}\log U,\\
			a_5&\coloneqq\frac1{p+1}\int_\R e^{-bt}U^{p+1}(t)\,dt>0.
		\end{aligned}
	\end{equation}

	Since \(e^{-t}U(t)\to\alpha_{N,s}\) as \(t\to-\infty\), define
	\begin{equation}\label{eq:a2}
		a_2\coloneqq\alpha_{N,s}\int_\R e^tU^{p_s^*}(t)\,dt>0.
	\end{equation}
	These integrals converge; in particular,
	\[
		e^tU^{p_s^*}(t)=
		\begin{cases}
			O(e^{({p_s^*}+1)t}),&t\to-\infty,\\
			O(e^{-({p_s^*}-1)t}),&t\to+\infty.
		\end{cases}
	\]
	The sign of \(a_4\) will not be used.

\subsection{Interaction asymptotics}

	\begin{lemma}\label{lem:interaction}
		As \(L\to\infty\),
		\begin{equation}\label{eq:interaction}
			\int_\R U^{p_s^*}(t)U(t-L)\,dt
			=a_2e^{-L}+o(e^{-L}).
		\end{equation}
	\end{lemma}

	\begin{proof}
		The bound \(U(z)\le Ce^z\) gives
		\[
			e^LU^{p_s^*}(t)U(t-L)\le Ce^tU^{p_s^*}(t)\in L^1(\R).
		\]
		For fixed \(t\),
		\[
			e^LU(t-L)\to\alpha_{N,s}e^t.
		\]
		Dominated convergence proves \eqref{eq:interaction}. The same estimate
		and a reflection give
		\[
			\int_\R(1+|t|)U^{p_s^*}(t)U(t\mp L)\,dt\le Ce^{-L}.\qedhere
		\]
	\end{proof}

	\begin{lemma}[Energy expansion]\label{lem:energy-expansion}
		Uniformly for \(\Lambda\in\mathcal K_\delta\),
		\[
			\E_\eps(V)
			=na_1+C_{n,\eps}+na_4\eps
			+\eps\Psi_n(\Lambda)+o(\eps),
		\]
		where
		\[
			C_{n,\eps}
			\coloneqq a_3\Big(\frac nb+\frac{n(n-1)}2\Big)\eps\log\eps
		\]
		and
		\[
			\Psi_n(\Lambda)
			=a_3n\log\Lambda_1-a_5\Lambda_1^b
			+\sum_{i=2}^n\bigl[a_3(n-i+1)\log\Lambda_i-a_2\Lambda_i\bigr].
		\]
	\end{lemma}

	\begin{proof}
		All constants are uniform on \(\mathcal K_\delta\).

		\noindent\emph{Step 1.} The critical part. Since \(\Ps U_i=U_i^{p_s^*}\),
		\[
			\frac12\Bform(V,V)=\frac12\sum_i\Bform(U_i,U_i)+\sum_{i<j}\int U_i^{p_s^*}U_j.
		\]
		On the cell \(\Omega_i\) of \(\xi_i\), set \(S_i\coloneqq\sum_{j\ne i}U_j\le CU_i\). Taylor's formula gives
		\[
			\Bigl|V^{{p_s^*}+1}-\sum_jU_j^{{p_s^*}+1}-({p_s^*}+1)\sum_{j\ne k}U_j^{p_s^*}U_k\Bigr|
			\le C\sum_{j\ne i}\bigl(U_i^{{p_s^*}-1}U_j^2+U_j^{p_s^*}U_i\bigr).
		\]
		For two centers at distance \(L\), the exponential tails give
		\[
			\int_{-\infty}^{L/2}\bigl(U^{{p_s^*}-1}U_L^2+U_L^{p_s^*}U\bigr)\le C(1+L)e^{-(1+\sigma_1)L},
		\]
		where \(U_L\coloneqq U(\cdot-L)\) and \(\sigma_1\coloneqq\min\{1,\frac{{p_s^*}-1}2\}\).
		Using \(|\xi_i-\xi_j|\ge|\log\eps|-C\) and \(\int U_i^{p_s^*}U_j=\Bform(U_i,U_j)=\int U_j^{p_s^*}U_i\), we obtain
		\[
			\frac12\Bform(V,V)-\frac1{{p_s^*}+1}\int V^{{p_s^*}+1}
			=na_1-\sum_{i<j}\int U_i^{p_s^*}U_j+O(\eps^{1+\sigma_1}|\log\eps|).
		\]
		By Lemma~\ref{lem:interaction} and \eqref{eq:int-balance}, adjacent pairs contribute \(-a_2\eps\Lambda_{i+1}+o(\eps)\), uniformly since \(\xi_{i+1}-\xi_i=|\log\eps|+O_\delta(1)\); for \(j\ge i+2\), \(e^{-(\xi_j-\xi_i)}=O(\eps^2)\).

		\noindent\emph{Step 2.} The supercritical perturbation. Set \(\ell_U(z)\coloneqq z+\log U(z)\). Then
		\[
			\ell_U\le C,\qquad |\ell_U(z)|\le C(1+|z|),\qquad
			|e^{\eps\ell_U}-1-\eps\ell_U|\le C\eps^2(1+|z|)^2.
		\]
		For \(t=z+\xi_i\),
		\[
			e^{\eps t}U_i^{{p_s^*}+1+\eps}(t)=e^{\eps\xi_i}e^{\eps\ell_U(z)}U^{{p_s^*}+1}(z).
		\]
		Taylor expansion, \(\eps\xi_i=O(\eps|\log\eps|)\), and \(\int_\R zU^{{p_s^*}+1}(z)\,dz=0\) give
		\[
			-\frac1{{p_s^*}+1+\eps}\int_\R e^{\eps t}U_i^{{p_s^*}+1+\eps}
			=-\frac1{{p_s^*}+1}\int_\R U^{{p_s^*}+1}-a_3\eps\xi_i+a_4\eps+O\bigl(\eps^2(1+\xi_i^2)\bigr).
		\]
		The cell expansion of Step~1 applies because
		\[
			e^{\eps t}U_i^\eps(t)=e^{\eps(t+\log U_i(t))}\le C,
			\qquad
			|e^{\eps t}U_i^{{p_s^*}+\eps}-U_i^{p_s^*}|
			\le C\eps(1+|\xi_i|+|z|)U_i^{p_s^*}.
		\]
		By the last bound in Lemma~\ref{lem:interaction}, the interaction terms change by \(O(\eps^2|\log\eps|)\). Thus the supercritical contribution is \(-a_3\eps\sum_i\xi_i+na_4\eps+o(\eps)\). Using \eqref{eq:xi-explicit},
		\[
			\sum_{i=1}^n\xi_i=-\Bigl(\frac nb+\frac{n(n-1)}2\Bigr)\log\eps-n\log\Lambda_1-\sum_{j=2}^n(n-j+1)\log\Lambda_j,
		\]
		which gives \(C_{n,\eps}\) and the logarithmic terms of \(\Psi_n\).

		\noindent\emph{Step 3.} The subcritical term. Translation gives
		\[
			\begin{gathered}
			-\frac1{p+1}\int_\R e^{-bt}U_i^{p+1}=-a_5e^{-b\xi_i},\\
			e^{-b\xi_1}=\eps\Lambda_1^b,\qquad e^{-b\xi_i}=O(\eps^{1+(i-1)b})\quad(i\ge2).
			\end{gathered}
		\]
		For the cross terms,
		\[
			0\le V^{p+1}-\sum_iU_i^{p+1}\le C\sum_{i\ne j}U_i^pU_j,\qquad
			U(z\mp L)\le Ce^{-L}e^{\pm z}.
		\]
		Thus
		\[
			\int_\R e^{-bt}U_i^pU_j\le Ce^{-b\xi_i}e^{-L}\int_\R e^{(\pm1-b)z}U^p(z)\,dz\le C\eps^2 ,
		\]
		since the last integrals converge: their slowest decay is \(e^{(p-b-1)z}\) as \(z\to-\infty\), with \(p-b-1=2p-{p_s^*}-1>0\) by \eqref{eq:2pP}. Thus
		\begin{equation}\label{eq:subcrit-energy-total}
			-\frac1{p+1}\int_\R e^{-bt}V^{p+1}
			=-a_5\eps\Lambda_1^b
			+O\bigl(\eps^{1+\min\{1,b\}}\bigr).
		\end{equation}
		Combining the three steps proves the lemma.
	\end{proof}

\subsection{The correction and the finite-dimensional critical point}

	Let \(\phi(\xi)\) be the solution of Proposition~\ref{prop:nonlinear}. Define
	\[
		J_\eps(\xi)\coloneqq\E_\eps(V+\phi(\xi)).
	\]

	\begin{lemma}\label{lem:J-E}
		Uniformly for \(\Lambda\in\mathcal K_\delta\),
		\[
			J_\eps(\xi(\Lambda))
			=\E_\eps(V)+O(\eps^{2\alpha})
			=\E_\eps(V)+o(\eps).
		\]
	\end{lemma}

	\begin{proof}
		Since \(\Ps V=\sum_iU_i^{p_s^*}\) and \(V>0\),
		\[
			D\E_\eps(V)[\psi]=-\int\mathfrak R_{\eps,\xi}\psi.
		\]
		The bound \(|\phi(t_\pm(y))|\le C\eps^\alpha(1+y)e^{-y}\) gives, on the left and right half-lines respectively,
		\[
			\begin{aligned}
				e^{-bt}(V+|\phi|)^{p-1}&\le C\eps(1+y)^{p-1}e^{-(2p-{p_s^*}-1)y},\\
				e^{\eps t}(V+|\phi|)^{{p_s^*}+\eps-1}&\le C(1+y)^{{p_s^*}+\eps-1}e^{-({p_s^*}-1)y}.
			\end{aligned}
		\]
		Both weights are bounded. Thus, for \(0\le\lambda\le1\),
		\[
			|D^2\E_\eps(V+\lambda\phi)[\phi,\phi]|\le\Bform(\phi,\phi)+C\|\phi\|_{L^2}^2.
		\]
		Testing \eqref{eq:projected-nl} with \(\phi\) gives \(\Bform(\phi,\phi)\le C\eps^{2\alpha}\). Also,
		\[
			\left|\int\mathfrak R\phi\right|+\|\phi\|_{L^2}^2
			\le\bigl(\|\mathfrak R\|_\gamma+\|\phi\|_\gamma\bigr)\|\phi\|_\gamma\int_\R\Wgt_\gamma^2
			\le C\eps^{2\alpha}.
		\]
		Taylor's formula gives \(|\E_\eps(V+\phi)-\E_\eps(V)|\le C\eps^{2\alpha}=o(\eps)\), since \(2\alpha>1\).
	\end{proof}

	\begin{lemma}[Variational reduction]\label{lem:critical-equivalence}
		For \(\eps\) sufficiently small, \(J_\eps\) is \(C^1\) on an open
		neighborhood of \(\{\xi(\Lambda):\Lambda\in\mathcal K_\delta\}\). For \(\xi=\xi(\Lambda)\), \(\Lambda\in\mathcal K_\delta\), the function
		\(V+\phi(\xi)\) is a positive solution of \eqref{eq:transformed} if and only if
		\(\nabla_\xi J_\eps(\xi)=0\).
	\end{lemma}

	\begin{proof}
		By Proposition~\ref{prop:nonlinear}, \(\xi\mapsto\phi(\xi)\) is \(C^1\) into \(\mathcal Y=H^s_k\cap X\). Proposition~\ref{prop:H-equivalence}(iii) and \(|U'''|\le Ce^{-|t|}\) give the same regularity for \(\xi\mapsto V\), with \(\partial_{\xi_\ell}V=-Z_\ell\). Since \(\E_\eps\in C^2(\mathcal Y)\), \(J_\eps\in C^1\).
		Comparing the projected equation with \eqref{eq:energy-derivatives} gives
		\[
			D\E_\eps(V+\phi)[\psi]=\sum_ic_i\int Z_i\psi,\qquad \psi\in\mathcal Y,
		\]
		the weak form of that equation. Taking \(\psi=\partial_{\xi_\ell}(V+\phi)\) and differentiating the orthogonality conditions, we obtain \(\partial_{\xi_\ell}J_\eps=\sum_ic_i\mathbf M_{i\ell}\), where
		\[
			\mathbf M_{i\ell}=-\int Z_iZ_\ell+\delta_{i\ell}\int\phi U_i''=-c_Z\delta_{i\ell}+o(1).
		\]
		Here \(\int\phi U_i''=O(\eps^\alpha)\) and \eqref{eq:Gram} give the last equality. Thus \(\mathbf M\) is invertible, and \(\nabla J_\eps=0\) if and only if \(c=0\). For \(c=0\), Lemma~\ref{lem:positive} gives \(V+\phi>0\), so the truncation can be removed. Conversely, if \(V+\phi\) is a positive solution of \eqref{eq:transformed}, then \(\sum_ic_iZ_i=0\), hence \(c=0\) by \eqref{eq:Gram}.
	\end{proof}

	\begin{lemma}[The limiting reduced energy]\label{lem:Psi}
		The function \(\Psi_n\) has a unique critical point in \((0,\infty)^n\):
		\begin{equation}\label{eq:Lambda-star}
			\Lambda_1^*=\Big(\frac{na_3}{ba_5}\Big)^{1/b},
			\qquad
			\Lambda_i^*=\frac{a_3(n-i+1)}{a_2},\quad 2\le i\le n.
		\end{equation}
		This point is a strict global maximum. Since \(\Psi_n\) is a sum of functions of one variable, its Hessian is diagonal, with
		\[
				\partial_{11}^2\Psi_n(\Lambda^*)
				=-\frac{a_3nb}{(\Lambda_1^*)^2}<0,\qquad
				\partial_{ii}^2\Psi_n(\Lambda^*)
				=-\frac{a_3(n-i+1)}{(\Lambda_i^*)^2}<0,
				\qquad i\ge2.
		\]
	\end{lemma}

	\begin{proof}
		Write \(\Psi_n(\Lambda)=\sum_i f_i(\Lambda_i)\), where
		\[
			f_1(x)=a_3n\log x-a_5x^b,
			\qquad
			f_i(x)=a_3(n-i+1)\log x-a_2x\quad(i\ge2).
		\]
		Then
		\[
			xf_1'(x)=a_3n-a_5bx^b,
			\qquad
			xf_i'(x)=a_3(n-i+1)-a_2x.
		\]
		Each \(xf_i'(x)\) decreases strictly and vanishes only at the value in \eqref{eq:Lambda-star}. Thus \(f_i\) increases up to that point and decreases thereafter. Also,
		\[
			\lim_{x\downarrow0}f_i(x)
			=\lim_{x\to\infty}f_i(x)=-\infty.
		\]
		Thus \(\Lambda^*\) is the unique critical point and the strict global maximum. This also covers \(0<b<1\), without global concavity of \(f_1\).
		Finally, \(a_5b(\Lambda_1^*)^b=a_3n\) gives
		\[
			f_1''(\Lambda_1^*)
			=-\frac{a_3n}{(\Lambda_1^*)^2}
			-a_5b(b-1)(\Lambda_1^*)^{b-2}
			=-\frac{a_3nb}{(\Lambda_1^*)^2},\qquad f_i''(\Lambda_i^*)
			=-\frac{a_3(n-i+1)}{(\Lambda_i^*)^2}.\qedhere
		\]
	\end{proof}

	\begin{remark}[Explicit constants]\label{rem:explicit}
		The substitution \(x=e^{2y}\) gives, for \(|\lambda|<m\),
		\[
			\int_\R\frac{e^{\lambda y}}{(2\cosh y)^{m}}\,dy
				=\frac12\int_0^\infty
				x^{(m+\lambda)/2-1}(1+x)^{-m}\,dx
				=\frac12B\left(\frac{m+\lambda}{2},
				\frac{m-\lambda}{2}\right).
		\]
		Using \(t=ay\) and \(U(ay)=\alpha_{N,s}(2\cosh y)^{-a}\), set
		\[
			\nu_p\coloneqq\frac{(N-2s)(2p+1-{p_s^*})}4>0.
		\]
		Since \(a({p_s^*}+1)=N\),
		\[
			\begin{aligned}
				\int_\R U^{{p_s^*}+1}\, dt
				=\frac a2\alpha_{N,s}^{{p_s^*}+1}B &\left(\frac N2,\frac N2\right),\quad
				\int_\R e^tU^{p_s^*}\,dt
				=\frac a2\alpha_{N,s}^{p_s^*} B\left(\frac N2,s\right),\\
				\int_\R &e^{-bt}U^{p+1}\,dt
				=\frac a2\alpha_{N,s}^{p+1}B\left(\nu_p,\frac N2\right).
			\end{aligned}
		\]
		Thus,
		\[
			\begin{aligned}
				&a_1=\frac{sa}{2N}\alpha_{N,s}^{{p_s^*}+1}B\left(\frac N2,\frac N2\right),\qquad
				&&a_2=\frac a2\alpha_{N,s}^{{p_s^*}+1}B\left(\frac N2,s\right),\\
				&a_3=\frac{a}{2({p_s^*}+1)}\alpha_{N,s}^{{p_s^*}+1}B\left(\frac N2,\frac N2\right),\qquad
				&&a_5=\frac{a}{2(p+1)}\alpha_{N,s}^{p+1}B\left(\nu_p,\frac N2\right).
			\end{aligned}
		\]
		Substitution in \eqref{eq:Lambda-star} gives
		\[
			\begin{aligned}
				\Lambda_i^*
				&=\frac{(N-2s)(n-i+1)}{2N}
				\frac{B\bigl(\frac N2,\frac N2\bigr)}
				{B\bigl(\frac N2,s\bigr)},
				\qquad i\ge2,\\[1mm]
				\Lambda_1^*
				&=\alpha_{N,s}\bigg(
				\frac{n(p+1)}{b({p_s^*}+1)}
				\frac{B\bigl(\frac N2,\frac N2\bigr)}
				{B\bigl(\nu_p,\frac N2\bigr)}\bigg)^{1/b}.
			\end{aligned}
		\]
		The amplitude \(\alpha_{N,s}\) is fixed by \eqref{eq:alpha-value}.
		By \eqref{eq:alpha-value} and Lemma~\ref{lem:indicial}(iii),
		\[
			a_2=\alpha_{N,s}^2\beta_s .
		\]
		The Green representation gives the same identity, as in the proof of Corollary~\ref{cor:asymptotics}(v). By Lemma~\ref{lem:indicial}(iv), \(U=G_s*U^{p_s^*}\), since \(U,U^{p_s^*}\in L^2\). Using \(e^{-t}U(t)\to\alpha_{N,s}\) as \(t\to-\infty\) and the evenness of \(U\),
		\[
			\alpha_{N,s}=\beta_s^{-1}\int_\R e^{-r}U^{p_s^*}(r)\,dr
			=\beta_s^{-1}\int_\R e^rU^{p_s^*}(r)\,dr
			=\frac{a_2}{\beta_s\alpha_{N,s}}.
		\]
	\end{remark}

\section{Proof of Theorem~\ref{thm:main}}\label{sec:proof-main}

	Fix \(n\ge1\) and the exponents in \eqref{eq:gamma-alpha-thm} and \eqref{eq:gamma-eta}. Choose \(\delta>0\) so that the point \(\Lambda^*\) in
	Lemma~\ref{lem:Psi} belongs to \(\operatorname{int}\mathcal K_\delta\).
	All constants below may depend on \(N,s,p,n,\gamma\) and \(\alpha\) (and on \(\delta\) and \(\eta\), which are fixed in terms of \(N,s,p,n\)), but not on \(\eps\).

	\begin{proof}[Proof of Theorem~\ref{thm:main}]
		\noindent\emph{Step 1.} Selection of the parameters.
		Take \(r>0\) such that
		\(\mathcal U\coloneqq\overline{B_r(\Lambda^*)}\subset\operatorname{int}\mathcal K_\delta\),
		and set
		\[
			\widetilde\Gamma_\eps(\Lambda)
			\coloneqq\frac{J_\eps(\xi(\Lambda))-na_1-C_{n,\eps}-na_4\eps}{\eps}.
		\]
		Lemmas~\ref{lem:energy-expansion} and \ref{lem:J-E} give
		\begin{equation}\label{eq:C0conv}
			\sup_{\Lambda\in\mathcal K_\delta}
			\bigl|\widetilde\Gamma_\eps(\Lambda)-\Psi_n(\Lambda)\bigr|\to0.
		\end{equation}
		The function \(\widetilde\Gamma_\eps\) is \(C^1\) on an open neighborhood of
		\(\mathcal U\). Since \(\Lambda^*\) is a strict maximum, the number
		\[
			\kappa_0\coloneqq\Psi_n(\Lambda^*)-\max_{\partial\mathcal U}\Psi_n
		\]
		is positive. If the supremum in \eqref{eq:C0conv} is smaller than \(\kappa_0/3\),
		then
		\[
			\max_{\partial\mathcal U}\widetilde\Gamma_\eps
			\le\Psi_n(\Lambda^*)-\frac{2\kappa_0}{3}
			<\widetilde\Gamma_\eps(\Lambda^*).
		\]
		Thus \(\widetilde\Gamma_\eps\) attains its maximum on \(\mathcal U\) at an interior
		point \(\Lambda_\eps\). In particular,
		\[
			\nabla_\Lambda\widetilde\Gamma_\eps(\Lambda_\eps)=0.
		\]
		Set \(\xi_\eps\coloneqq\xi(\Lambda_\eps)\). The chain rule gives
		\[
			0=\nabla_\Lambda\widetilde\Gamma_\eps(\Lambda_\eps)
			=\eps^{-1}\left(\frac{\partial\xi}{\partial\Lambda}(\Lambda_\eps)\right)^T\nabla_\xi J_\eps(\xi_\eps).
		\]
		The Jacobian is invertible, so \(\nabla_\xi J_\eps(\xi_\eps)=0\). By Lemma~\ref{lem:critical-equivalence}, \(c_1=\cdots=c_n=0\).

		By compactness of \(\mathcal U\), it suffices to consider any convergent subsequence \(\Lambda_{\eps_j}\to\overline\Lambda\in\mathcal U\). Uniform convergence in \eqref{eq:C0conv} and continuity of \(\Psi_n\) give
		\[
			\Psi_n(\overline\Lambda)
			=\lim_j\widetilde\Gamma_{\eps_j}(\Lambda_{\eps_j})
			\ge\lim_j\widetilde\Gamma_{\eps_j}(\Lambda^*)
			=\Psi_n(\Lambda^*).
		\]
		The maximum is unique, so \(\overline\Lambda=\Lambda^*\). Thus \(\Lambda_\eps\to\Lambda^*\).

		\medskip

		\noindent\emph{Step 2.} Return to the original variables.
		Let
		\[
			\phi_{\eps,n}\coloneqq\phi(\xi_\eps),\qquad
			\mu_{i,\eps}\coloneqq e^{-\xi_{i,\eps}/a},
		\]
		and set
		\[
			v_{\eps,n}\coloneqq\sum_{i=1}^nU(\,\cdot-\xi_{i,\eps})+\phi_{\eps,n}.
		\]
		Proposition~\ref{prop:nonlinear} gives
		\begin{equation}\label{eq:remainder-main}
			\|\phi_{\eps,n}\|_\gamma\le C\eps^\alpha.
		\end{equation}
		By Lemma~\ref{lem:positive},
		\[
			u_{\eps,n}(x)
			\coloneqq|x|^{-a}v_{\eps,n}(-a\log|x|),\qquad x\ne0,
		\]
		extends to \(0\) as a positive radial classical solution of
		\eqref{eq:main-problem}. It is bounded and belongs to
		\(C^{2s+\sigma}(\R^N)\) for every \(0<\sigma<\min\{1,2-2s\}\) with \(2s+\sigma\ne1\).
		The decay estimate \eqref{eq:u-decay} is \eqref{eq:fast-decay}.
		By \eqref{eq:translation-scale},
		\[
			u_{\eps,n}(x)=\sum_{i=1}^nW_{\mu_{i,\eps},0}(x)
			+\omega_{\eps,n}(x),
			\qquad
			\omega_{\eps,n}(x)=e^t\phi_{\eps,n}(t),\quad t=-a\log|x|.
		\]
		For each \(i\), set \(R_i\coloneqq|x|/\mu_{i,\eps}\). Then
		\[
			t-\xi_{i,\eps}=-a\log R_i,
			\qquad
			e^te^{-\gamma|t-\xi_{i,\eps}|}
			=\mu_{i,\eps}^{-a}
			\min\{R_i^{-a(1-\gamma)},R_i^{-a(1+\gamma)}\}.
		\]
		Together with \eqref{eq:remainder-main}, this proves
		\eqref{eq:remainder-x}.

		\medskip

		\noindent\emph{Step 3.} Scales and energy.
		The definitions of the centers give
		\[
			e^{-\xi_{1,\eps}}=\eps^{1/b}\Lambda_{1,\eps},
			\qquad
			e^{-(\xi_{i,\eps}-\xi_{i-1,\eps})}=\eps\Lambda_{i,\eps}
			\quad(2\le i\le n).
		\]
		Raising these identities to the power \(1/a\) proves \eqref{eq:mu1} and \eqref{eq:muratio}.
		Multiplication gives
		\[
			\mu_{i,\eps}
			=\bigg(\prod_{j=1}^i\Lambda_{j,\eps}\bigg)^{1/a}
			\eps^{\frac1a(\frac1b+i-1)}.
		\]
		The convergence \(\Lambda_\eps\to\Lambda^*\) and Lemma~\ref{lem:Psi}
		give \eqref{eq:Lambda-star-intro} and \eqref{eq:mu-asymp-intro}.
		Finally, Lemmas~\ref{lem:energy-expansion} and~\ref{lem:J-E} give
		\[
			\E_\eps(v_{\eps,n})=na_1+O(\eps|\log\eps|).\qedhere
		\]

	\end{proof}

\subsection{Asymptotic consequences}

	\begin{proof}[Proof of Corollary~\ref{cor:asymptotics}]
		Throughout the proof, \(\mu_i=\mu_{i,\eps}\), \(\xi_i=\xi_{i,\eps}\), \(v=v_{\eps,n}\), \(V=\sum_iU(\cdot-\xi_{i,\eps})\) and \(q={p_s^*}+\eps\).

		\medskip
		\noindent (i) Multiply \eqref{eq:relative-error-EF} by \(e^t\) and use continuity at \(x=0\).

		\medskip
		\noindent (ii) By Lemma~\ref{lem:positive}(c), \(u_{\eps,n}\) and the bubbles are continuous at \(0\). Thus \(\omega_{\eps,n}(0)=\lim_{t\to\infty}e^t\phi_{\eps,n}(t)\) exists, and \eqref{eq:phi-tail-improved} gives
		\[
			|\omega_{\eps,n}(0)|\le C\eps^\alpha e^{\xi_n}=C\eps^\alpha\mu_n^{-a}.
		\]
		Together with \(W_{\mu,0}(0)=\alpha_{N,s}\mu^{-a}\), this proves (ii).

		\medskip
		\noindent (iii) For \(j\ne i\), the rescaled bubble is
		\[
			\mu_i^aW_{\mu_j,0}(\mu_i y)
			=\alpha_{N,s}\frac{(\mu_j/\mu_i)^a}{((\mu_j/\mu_i)^2+|y|^2)^a}.
		\]
		It tends to zero uniformly on compact subsets of \(\R^N\setminus\{0\}\), and of \(\R^N\) when \(i=n\). Part (i) controls the remainder.

		\medskip
		\noindent (iv) Polar coordinates give
		\[
			\begin{aligned}
				\int_{\R^N}u_{\eps,n}^{{p_s^*}+1}\,dx&=\frac{|\Sph^{N-1}|}a\int_\R v^{{p_s^*}+1}\,dt,\\
				\int_{\R^N}W_{1,0}^{{p_s^*}+1}\,dx&=\frac{|\Sph^{N-1}|}a\int_\R U^{{p_s^*}+1}\,dt.
			\end{aligned}
		\]
		By \eqref{eq:relative-error-EF},
		\[
			|v^{{p_s^*}+1}-V^{{p_s^*}+1}|\le C\eps^\varrho V^{{p_s^*}+1}.
		\]
		Step~1 of the proof of Lemma~\ref{lem:energy-expansion} and Lemma~\ref{lem:interaction} give
		\[
			\int_\R V^{{p_s^*}+1}=n\int_\R U^{{p_s^*}+1}+O(\eps),
		\]
		which proves (iv).

		\medskip
		\noindent (v) Set \(H=e^{-bt}v^p+e^{\eps t}v^q\). As in the proof of Lemma~\ref{lem:positive}(c), \(v=G_s*H\). On \(|t-r|\ge1\),
		\[
			e^{-t}G_s(t-r)\le Ce^{-r},\qquad
			e^{-t}G_s(t-r)\to\beta_s^{-1}e^{-r}\quad(t\to-\infty).
		\]
		Since \(e^{-r}H\in L^1\), dominated convergence applies. On \(|t-r|<1\),
		\[
			e^{-t}\int_{|t-r|<1}G_s(t-r)H(r)\,dr
			\le e\|G_s\|_{L^1}\sup_{r<t+1}e^{-r}H(r)\to0,
		\]
		because
		\[
			e^{-r}H(r)=O\bigl(e^{(2p-{p_s^*}-1)r}+e^{({p_s^*}+2\eps-1)r}\bigr)\qquad(r\to-\infty).
		\]
		Thus
		\[
			A_{\eps,n}=\lim_{t\to-\infty}e^{-t}v(t)=\beta_s^{-1}\int_\R e^{-r}H(r)\,dr>0.
		\]
		Comparing with \(\lim_{t\to-\infty}e^{-t}V(t)=\alpha_{N,s}\sum_i\mu_i^a\) and using \eqref{eq:phi-tail-improved} gives the expansion.
	\end{proof}

	\begin{proof}[Proof of Corollary~\ref{cor:multiplicity}]
		Fix \(m\ge1\) and construct the solutions for \(n=1,\ldots,m\).
		Let
		\[
			M_*\coloneqq\int_{\R^N}W_{1,0}^{{p_s^*}+1}\,dx>0.
		\]
		By Corollary~\ref{cor:asymptotics}(iv), choose \(\eps_m'\in(0,\min_{n\le m}\eps_n]\) such that, for \(0<\eps<\eps_m'\),
		\[
			\left|\int_{\R^N}u_{\eps,n}^{{p_s^*}+1}\,dx-nM_*\right|<\frac{M_*}{3}
			\qquad(1\le n\le m).
		\]
		The intervals \((nM_*-\frac{M_*}3,\,nM_*+\frac{M_*}3)\), \(1\le n\le m\), are pairwise disjoint. Thus the solutions are distinct for every \(\eps\in(0,\eps_m')\).
	\end{proof}

\subsection{The balance between the two powers}

	Translation invariance of the quadratic form gives a Pohozaev-type identity on the cylinder. For fractional Dirichlet problems on bounded domains, the identity includes a boundary term; see \cite{RosOtonSerraPohozaev2014}.

	\begin{lemma}\label{lem:v-C1}
		Let \(v=v_{\eps,n}\) and \(H\coloneqq e^{-bt}v^p+e^{\eps t}v^q\), with \(q={p_s^*}+\eps\). Then \(H\in L^1(\R)\cap L^2(\R)\) and \(v\in C^{1,\theta}(\R)\) for some \(\theta\in(0,1)\); in particular, \(v'\in L^\infty(\R)\).
	\end{lemma}

	\begin{proof}
		The integrability of \(H\) follows from \eqref{eq:v-tail}. Choose
		\[
			\max\{0,1-2s\}<\sigma<\min\{1,2-2s\}.
		\]
		By Proposition~\ref{prop:nonlinear}, \(\phi_{\eps,n}\in C^{2s+\sigma}(\R)\). Since \(V=\sum_iU(\cdot-\xi_{i,\eps})\) is smooth with bounded derivatives, \(v=V+\phi_{\eps,n}\in C^{1,\theta}(\R)\), with \(\theta=2s+\sigma-1\in(0,1)\). Thus \(v'\in L^\infty(\R)\).
	\end{proof}

	\begin{proposition}\label{prop:balance}
		Let \(u=u_{\eps,n}\) be a solution constructed above, \(v=v_{\eps,n}\) and
		\(q={p_s^*}+\eps\). Then
		\begin{equation}\label{eq:balance-cylinder}
			\frac{b}{p+1}\int_\R e^{-bt}v^{p+1}\,dt
			=\frac{\eps}{q+1}\int_\R e^{\eps t}v^{q+1}\,dt.
		\end{equation}
		Equivalently,
		\begin{equation}\label{eq:balance-original}
			\frac{b}{p+1}\int_{\R^N}u^{p+1}\,dx
			=\frac{\eps}{q+1}\int_{\R^N}u^{q+1}\,dx.
		\end{equation}
	\end{proposition}

	\begin{proof}
		By Lemma~\ref{lem:v-C1}, \(H\in L^1(\R)\cap L^2(\R)\), \(v\in C^1(\R)\) and \(v'\in L^\infty(\R)\). Since \(v\in H^s_k\) and \(\Ps v=H\), Proposition~\ref{prop:H-equivalence}(ii) and density give
		\[
			\Bform(v,\psi)=\int_\R H\psi,\qquad \psi\in H^s_k(\R).
		\]
		Set \(\tau_hv\coloneqq v(\cdot+h)\). Translation invariance and symmetry give
		\[
			\Bform(v,\tau_hv)=\Bform(\tau_{-h}v,v)=\Bform(v,\tau_{-h}v),
		\]
		hence
		\[
			\int_\R H(t)\frac{v(t+h)-v(t-h)}{2h}\,dt=0\qquad(h\ne0).
		\]
		The difference quotients converge to \(v'\) and are bounded by \(\|v'\|_{L^\infty}\). Since \(H\in L^1\), dominated convergence gives \(\int_\R Hv'=0\). By \eqref{eq:v-tail}, the boundary terms \(e^{-bt}v^{p+1}\) and \(e^{\eps t}v^{q+1}\) vanish at both ends. Integration by parts on \([-R,R]\), followed by \(R\to\infty\), yields
		\[
			0=\int_\R Hv'=\frac{b}{p+1}\int_\R e^{-bt}v^{p+1}\,dt-\frac{\eps}{q+1}\int_\R e^{\eps t}v^{q+1}\,dt.
		\]
		This is \eqref{eq:balance-cylinder}; polar coordinates give \eqref{eq:balance-original}.
	\end{proof}

	\begin{remark}
		The balance \eqref{eq:balance-cylinder} also determines the first limiting parameter. By \eqref{eq:relative-error-EF}, \eqref{eq:subcrit-energy-total} and Steps~1--2 of the proof of Lemma~\ref{lem:energy-expansion}, applied to \(V\),
		\[
			\begin{aligned}
				\frac{b}{p+1}\int_\R e^{-bt}v^{p+1}&=b\,a_5\eps\Lambda_{1,\eps}^b(1+o(1)),\\
				\frac{\eps}{q+1}\int_\R e^{\eps t}v^{q+1}&=\eps n a_3(1+o(1)).
			\end{aligned}
		\]
		Here \eqref{eq:relative-error-EF} gives \(|v^r-V^r|\le C\eps^\varrho V^r\) for \(r=p+1,q+1\). Thus \(b\,a_5\Lambda_{1,\eps}^b\to n\,a_3\), which gives \(\Lambda_1^*\) in \eqref{eq:Lambda-star}.
	\end{remark}

\appendix
\section{Function spaces, weighted estimates and regularity}\label{app:functional}

	All results in this appendix hold for \(0<s<1\).

	\begin{proposition}\label{prop:H-equivalence}
		The space \(H^s_k(\R)\) of Section~\ref{sec:linear} coincides
		with \(H^s(\R)\), and
		\[
			\Bform(v,v)=\frac1{2\pi}\int_\R\Psym(k)|\widehat v(k)|^2\,dk,
			\qquad
			c\|v\|_{H^s}^2\le\Bform(v,v)\le C\|v\|_{H^s}^2.
		\]
		The following properties hold:
		\begin{enumerate}[label=\textup{(\roman*)}]
			\item The embedding \(H^s_k(\R)\hookrightarrow L^2_{\rm loc}(\R)\) is compact.
			\item For \(v\in H^s_k(\R)\) and \(\varphi\in C_c^\infty(\R)\),
			\(\Bform(v,\varphi)=\int_\R v\Ps\varphi\).
			In particular, a weak equation \(\Bform(v,\psi)=\int g\psi\), with
			\(g\in L^2(\R)\), implies \(\Ps v=g\) in \(\mathcal D'(\R)\).
			\item \(U_i,Z_i,U_i''\in H^s_k(\R)\cap X\), and
			\(\Wgt_\beta\in L^2(\R)\) for every \(\beta>0\).
		\end{enumerate}
	\end{proposition}

	\begin{proof}
		By Tonelli's theorem, Plancherel and Lemma~\ref{lem:indicial},
		\[
			\Bform(v,v)=\frac1{2\pi}\int_\R\Psym(k)|\widehat v(k)|^2\,dk
			\asymp\int_\R(1+k^2)^s|\widehat v(k)|^2\,dk
			\asymp\|v\|_{H^s}^2.
		\]
		For (i), apply the fractional Rellich theorem \cite[Theorem~7.1]{DNPV} on bounded intervals, which are extension domains by \cite[Theorem~5.4]{DNPV}. Assertions (ii)--(iii) are standard.
	\end{proof}

	\begin{lemma}[Bounded representatives]\label{lem:bounded-rep}
		Let \(w\in L^2(\R)\), \(\mathcal V\in L^\infty(\R)\) and \(f_1\in L^2(\R)\cap L^\infty(\R)\) satisfy \(\Ps w=\mathcal Vw+f_1\) in \(\mathcal D'(\R)\). Then \(w=G_s*(\mathcal Vw+f_1)\), \(w\) has a bounded uniformly continuous representative, and
		\[
			\|w\|_{L^\infty}
			\le C\bigl(\|w\|_{L^2}+\|f_1\|_{L^2}+\|f_1\|_{L^\infty}\bigr),
		\]
		where \(C\) depends only on \(N\), \(s\) and \(\|\mathcal V\|_{L^\infty}\).
	\end{lemma}

	\begin{proof}
		Set \(g\coloneqq\mathcal Vw+f_1\in L^2\). By Lemma~\ref{lem:indicial}(iv), \(w=G_s*g\). Fix \(\varsigma\in(1,\min\{2,\varsigma_s\})\), let \(\varsigma'\) be its conjugate exponent, and set \(\kappa\coloneqq1-1/\varsigma\in(0,\frac12)\). Since \(G_s\in L^1\cap L^{\varsigma}\),
		\[
			G_s*f_1\in L^2\cap L^\infty,
			\qquad \|G_s*f_1\|_{L^\infty}\le\|G_s\|_{L^1}\|f_1\|_{L^\infty}.
		\]
		Start with \(r_0=2\). If \(w\in L^2\cap L^{r_j}\) and \(1/r_j>\kappa\), Young's inequality gives
		\[
			G_s*(\mathcal Vw)\in L^{r_{j+1}},
			\qquad \frac1{r_{j+1}}\coloneqq\frac1{r_j}-\kappa.
		\]
		Thus \(w\in L^2\cap L^{r_{j+1}}\). After finitely many steps, \(r_j\ge\varsigma'\ge2\). Interpolation and H\"older's inequality then give
		\[
			\|G_s*(\mathcal Vw)\|_{L^\infty}\le\|G_s\|_{L^{\varsigma}}\|\mathcal Vw\|_{L^{\varsigma'}}.
		\]
		Translation continuity in \(L^1\) and \(L^{\varsigma}\) gives uniform continuity of both convolutions.
	\end{proof}

	\begin{lemma}[Weighted resolvent estimates]\label{lem:resolvent}
		Let \(0<\beta<1\), \(w\in L^2(\R)\), and \(\Ps w=f\) in \(\mathcal D'(\R)\), with \(\|f\|_\beta<\infty\). Then
		\[
			\|w\|_\beta\le C\|f\|_\beta,
		\]
		where \(C\) depends only on \(N\), \(s\) and \(\beta\), and is independent of the centers. The condition \(w\in L^2(\R)\) excludes the homogeneous solutions \(e^{\pm t}\).
	\end{lemma}

	\begin{proof}
		Since \(|f|\le\|f\|_\beta\Wgt_\beta\) and \(\Wgt_\beta\in L^2\), \(f\in L^2\). Lemma~\ref{lem:indicial}(iv) gives \(w=G_s*f\). The weight \(\Wgt_\beta\) is slowly varying. By Lemma~\ref{lem:comparison},
		\[
			|w(t)|\le C\|f\|_\beta\int_\R e^{-|t-r|}\Wgt_\beta(r)\,dr
			\le C\|f\|_\beta\Wgt_\beta(t),
		\]
		where the last inequality follows from the first bound in \eqref{eq:conv-elementary}, applied to each summand.
	\end{proof}

	The bounds \eqref{eq:kernel-bounds} verify the kernel hypotheses of the Schauder estimate \cite[Proposition~3.9]{MR3694655}: comparability to \(|\tau|^{-1-2s}\) for \(|\tau|\le1\), and polynomial decay for \(|\tau|>1\). The proof of \cite[Proposition~3.9]{MR3694655} combines \cite[Theorem~1.2]{DongKim2013} with \cite[Proposition~2.8]{SilvestreCPAM2007}.

	\begin{remark}
		By \eqref{eq:kernel-expansion}, \(k_s=c_1|\tau|^{-1-2s}+\tilde k\) with \(c_1=a^{2s}c_{1,s}\) and \(\tilde k\in L^1(\R)\) even. Hence
		\[
			\Ps w=a^{2s}(-\Delta)^s_{\R}w+c_*w-\tilde k*w,
			\qquad c_*\coloneqq m_s+\int_\R\tilde k,
		\]
		where \((-\Delta)^s_\R\) is the one-dimensional fractional Laplacian. Since \(\|\tilde k*w\|_{C^{0,\sigma}(\R)}\le\|\tilde k\|_{L^1}\|w\|_{C^{0,\sigma}(\R)}\), Lemma~\ref{lem:schauder} also follows from \cite[Proposition~2.8]{SilvestreCPAM2007} in dimension one. Also,
		\[
			\Psym(k)=a^{2s}|k|^{2s}+c_*-\widehat{\tilde k}(k),
			\qquad \bigl|\Psym(k)-a^{2s}|k|^{2s}\bigr|\le|c_*|+\|\tilde k\|_{L^1}.
		\]
		Together with \(\Psym\ge m_s\), this gives \(c(1+|k|^{2s})\le\Psym(k)\le C(1+|k|^{2s})\), as in Lemma~\ref{lem:indicial}(ii).
	\end{remark}

	\begin{lemma}[Schauder estimate on the line]\label{lem:schauder}
		Let \(0<\sigma<\min\{1,2-2s\}\) with \(2s+\sigma\ne1\), and let \(w,f\in C^{0,\sigma}(\R)\) satisfy \(\Ps w=f\) in \(\mathcal D'(\R)\). Then \(w\in C^{2s+\sigma}(\R)\) 
		and
		\[
			\|w\|_{C^{2s+\sigma}(\R)}\le C\left(\|w\|_{C^{0,\sigma}(\R)}+\|f\|_{C^{0,\sigma}(\R)}\right),
		\]
		where \(C\) depends only on \(N\), \(s\) and \(\sigma\).
	\end{lemma}

	\begin{proof}
		For smooth \(w\) with bounded derivatives, the equation holds pointwise and \cite[Proposition~3.9]{MR3694655} applies. For general \(w\), set \(w_\delta\coloneqq\rho_\delta*w\), with \(\rho_\delta\) an even mollifier. Translation invariance gives, for \(\varphi\in C_c^\infty(\R)\),
		\[
			\int_\R w_\delta\Ps\varphi=\int_\R w\Ps(\rho_\delta*\varphi)=\int_\R f\,(\rho_\delta*\varphi)=\int_\R(\rho_\delta*f)\varphi.
		\]
		By the truncation argument in Lemma~\ref{lem:indicial}(iii) (with truncation parameter \(\tau\to0\) for fixed \(\delta\)), \(\int(\Ps w_\delta)\varphi=\int w_\delta\Ps\varphi\). Both \(\Ps w_\delta\) and \(\rho_\delta*f\) are continuous, so \(\Ps w_\delta=\rho_\delta*f\) pointwise. Moreover,
		\[
			\|w_\delta\|_{C^{0,\sigma}(\R)}\le\|w\|_{C^{0,\sigma}(\R)},
			\qquad \|\rho_\delta*f\|_{C^{0,\sigma}(\R)}\le\|f\|_{C^{0,\sigma}(\R)}.
		\]
		The smooth estimate is therefore uniform in \(\delta\). Since \(w_\delta\to w\) locally uniformly, Arzel\`a--Ascoli and lower semicontinuity of the norm give the result.
	\end{proof}

	\begin{proposition}[Regularity on the line]\label{prop:regularity}
		Let \(w,f\in L^\infty(\R)\) satisfy \(\Ps w=f\) in \(\mathcal D'(\R)\).
		\begin{enumerate}[label=\textup{(\roman*)}]
			\item \(w=G_s*f\), and for every \(\beta\in(0,\min\{1,2s\})\),
			\[
				\|w\|_{C^{0,\beta}(\R)}\le C\|f\|_{L^\infty(\R)}.
			\]
			\item If, in addition, \(f\in C^{0,\sigma}(\R)\) with \(0<\sigma<\min\{1,2-2s\}\) and \(2s+\sigma\ne1\), then \(w\in C^{2s+\sigma}(\R)\),
			\[
				\|w\|_{C^{2s+\sigma}(\R)}\le C\|f\|_{C^{0,\sigma}(\R)},
			\]
			and \(\Ps w=f\) pointwise in \(\R\).
		\end{enumerate}
		The constants depend only on \(N\), \(s\) and the H\"older exponents.
	\end{proposition}

	\begin{proof}
		\noindent\emph{Step 1.} Green representation. By the evenness of \(G_s\) and Lemma~\ref{lem:indicial}(ii)--(iii), for \(\varphi\in C_c^\infty(\R)\),
		\[
			\int_\R(G_s*f)\Ps\varphi=\int_\R f\,(G_s*\Ps\varphi)=\int_\R f\varphi.
		\]
		Also, \(\|G_s*f\|_{L^\infty}\le\|G_s\|_{L^1}\|f\|_{L^\infty}\). Thus \(w-G_s*f\) is bounded and satisfies \(\Ps(w-G_s*f)=0\). Lemma~\ref{lem:indicial}(iv) gives \(w=G_s*f\).

		\medskip

		\noindent\emph{Step 2.} H\"older continuity. Fix \(\beta\in(0,\min\{1,2s\})\). The \(L^\infty\) bound gives the estimate for \(|h|>\frac14\). For \(0<|h|\le\frac14\), Lemma~\ref{lem:indicial}(iii) gives
		\[
			\int_{|z|\le2|h|}\bigl(G_s(z+h)+G_s(z)\bigr)\,dz\le2\int_{|z|\le3|h|}G_s\le C\int_0^{3|h|}\varpi_s\le C_\beta|h|^\beta.
		\]
		For \(|z|>2|h|\), the segment from \(z\) to \(z+h\) lies in \(\{|\zeta|\ge|z|/2\}\). The mean value theorem and the bounds on \(G_s'\) give
		\[
			|G_s(z+h)-G_s(z)|\le C|h|\big(|z|^{2s-2}\mathbf 1_{\{|z|\le2\}}+e^{-|z|/2}\big),
		\]
		and \(\int_{2|h|<|z|\le2}|z|^{2s-2}\,dz\le C_\beta|h|^{\beta-1}\). Thus \(\|G_s(\cdot+h)-G_s\|_{L^1}\le C_\beta|h|^\beta\), and
		\[
			|w(t+h)-w(t)|=\bigg|\int_\R\bigl(G_s(t+h-r)-G_s(t-r)\bigr)f(r)\,dr\bigg|\le C_\beta|h|^\beta\|f\|_{L^\infty}.
		\]
		Together with Step~1, this proves (i).

		\medskip

		\noindent\emph{Step 3.} Higher regularity. Choose \(\beta_0\in(0,\min\{1,2s\})\) with \(\beta_0+2sj\ne1\) for every integer \(j\ge0\). By (i), 
		\(w\in C^{0,\beta_0}(\R)\). Suppose that \(w\in C^{0,\beta_j}(\R)\) with \(\beta_j\le1\), and set \(\sigma_j\coloneqq\min\{\sigma,\beta_j\}\). Then 
		\(w,f\in C^{0,\sigma_j}(\R)\) and \(2s+\sigma_j\ne1\), since either \(\sigma_j=\sigma\) or \(\sigma_j=\beta_j=\beta_0+2sj\). Lemma~\ref{lem:schauder} gives
		\[
			\|w\|_{C^{2s+\sigma_j}(\R)}\le C\big(\|w\|_{C^{0,\beta_j}(\R)}+\|f\|_{C^{0,\sigma}(\R)}\big).
		\]
		If \(\sigma_j=\sigma\), this gives (ii). Otherwise, set \(\beta_{j+1}\coloneqq\min\{1,\beta_j+2s\}\). When \(\beta_j+2s>1\),
		\[
			w\in C^{1,\beta_j+2s-1}(\R)\subset C^{0,1}(\R).
		\]
		The exponent increases by \(2s\) until it reaches \(1>\sigma\). Thus \(\sigma_j=\sigma\) after finitely many steps, and (i) gives the bound in (ii).

		Since \(2s+\sigma<2\),
		\[
			|2w(t)-w(t+z)-w(t-z)|\le C|z|^{2s+\sigma},
			\qquad t\in\R,\quad |z|\le1.
		\]
		Thus the principal value defining \(\Ps w\) converges locally uniformly. The truncation identity in Lemma~\ref{lem:indicial}(iii) identifies it with the distributional operator. Since both sides are continuous, \(\Ps w=f\) pointwise.
	\end{proof}

	\begin{lemma}[H\"older transfer]\label{lem:holder-transfer}
		Let \(r>1\), let \(\mathfrak w\in\{e^{-bt},e^{\eps t}\}\) with \(0<\eps\le1\), and let \(\mathcal M\colon\R\to(0,\infty)\) satisfy \(\mathcal M(t')\le K_0\mathcal M(t)\) whenever \(|t-t'|\le1\), and \(\mathfrak w\mathcal M^{r-1}+\mathcal M\le K_0\) on \(\R\). If \(w\in C^{0,\beta}(\R)\), with \(0<\beta\le1\), and \(|w|\le\mathcal M\), then
		\[
			\|\mathfrak w\,w_+^r\|_{C^{0,\beta}(\R)}\le C\bigl(1+[w]_{C^{0,\beta}(\R)}\bigr),
		\]
		where \(C\) depends only on \(r\), \(K_0\) and \(b\).
	\end{lemma}

	\begin{proof}
		The bound \(\mathfrak w w_+^r\le\mathfrak w\mathcal M^{r-1}\mathcal M\le K_0^2\) also controls differences with \(|t-t'|\ge1\). For \(|t-t'|<1\),
		\begin{align*}
			\mathfrak w(t')&\le C\mathfrak w(t),&
			|\mathfrak w(t)-\mathfrak w(t')|&\le C\mathfrak w(t)|t-t'|,\\
			|w_+^r(t)-w_+^r(t')|&\le rK_0^{r-1}\mathcal M(t)^{r-1}|w(t)-w(t')|.
		\end{align*}
		Thus
		\begin{align*}
			|\mathfrak w w_+^r(t)-\mathfrak w w_+^r(t')|
			&\le C\mathfrak w(t)\mathcal M(t)^r|t-t'|
			+C\mathfrak w(t)\mathcal M(t)^{r-1}[w]_{C^{0,\beta}}|t-t'|^\beta\\
			&\le C(1+[w]_{C^{0,\beta}})|t-t'|^\beta,
		\end{align*}
		which gives the stated norm bound.
	\end{proof}

	The change of variables in Lemma~\ref{lem:indicial}(iii) gives, for measurable \(h\) with \(G_s*|h|<\infty\) and \(x\ne0\),
	\begin{equation}\label{eq:Riesz-cylinder}
		\begin{aligned}
		|x|^{-a}(G_s*h)(-a\log|x|)
		&=c^{\mathrm R}_{N,s}\int_{\R^N}\frac{F(y)}{|x-y|^{N-2s}}\,dy,\\
		F(y)&\coloneqq|y|^{-a-2s}h(-a\log|y|).
		\end{aligned}
	\end{equation}
	Thus \(v=G_s*h\) corresponds to the Riesz potential of \(F\). For nonnegative distributional solutions of the fractional H\'enon equation, see \cite[Lemma~2.1]{BarriosQuaas2020} for a different proof.

	\begin{lemma}[Riesz potentials]\label{lem:Riesz}
		Let \(f\in L^1(\R^N)\cap L^\infty(\R^N)\) and \(u\coloneqq c^{\mathrm R}_{N,s}|\cdot|^{2s-N}*f\). Then, for every \(\beta\in(0,\min\{1,2s\})\),
		\[
			\|u\|_{C^{0,\beta}(\R^N)}\le C\left(\|f\|_{L^1(\R^N)}+\|f\|_{L^\infty(\R^N)}\right),
		\]
		where \(C\) depends only on \(N\), \(s\) and \(\beta\).
	\end{lemma}

	\begin{proof}
		Set \(\Phi(z)\coloneqq|z|^{2s-N}\). Splitting at \(|x-y|=1\) gives \(\|u\|_{L^\infty}\le C(\|f\|_{L^1}+\|f\|_{L^\infty})\), which gives the estimate for \(|h|>\frac14\). For \(0<|h|\le\frac14\),
		\[
			|u(x+h)-u(x)|\le C\|f\|_{L^\infty}\int_{|z|\le2}|\Phi(z+h)-\Phi(z)|+C\|f\|_{L^1}\sup_{|z|>2}|\Phi(z+h)-\Phi(z)|.
		\]
		Since \(|\nabla \Phi(\zeta)|\le C|\zeta|^{2s-N-1}\), the last supremum is at most \(C|h|\). As in Step~2 of the proof of 
		Proposition~\ref{prop:regularity},
		\begin{align*}
			\int_{|z|\le2|h|}\bigl(\Phi(z+h)+\Phi(z)\bigr)\,dz&\le C|h|^{2s},\\
			\int_{2|h|<|z|\le2}|\Phi(z+h)-\Phi(z)|\,dz&\le C|h|\int_{2|h|}^2r^{2s-2}\,dr\le C_\beta|h|^\beta,
		\end{align*}
		which proves the lemma.
	\end{proof}

	Let \(\phi\) be smooth and bounded with bounded derivatives, and set \(g\coloneqq(-\Delta)^s\phi\). If \(g\in C^{0,\sigma}(\R^N)\), with \(0<\sigma<\min\{1,2-2s\}\) and \(2s+\sigma\ne1\), \cite[Proposition~2.8]{SilvestreCPAM2007} gives
	\begin{equation}\label{eq:Silvestre-holder}
		\|\phi\|_{C^{2s+\sigma}(\R^N)}
		\le C\bigl(\|\phi\|_{L^\infty(\R^N)}+\|g\|_{C^{0,\sigma}(\R^N)}\bigr),
	\end{equation}
	where \(C\) depends only on \(N\), \(s\) and \(\sigma\). See also \cite[(3.9)--(3.10)]{MR3694655} for the one-dimensional case.
	For \(\phi\in L^\infty(\R^N)\) with \((-\Delta)^s\phi=g\) in \(\mathcal D'(\R^N)\), set \(\phi_\delta\coloneqq\rho_\delta*\phi\). These functions are smooth and bounded with all their derivatives. The duality argument in Lemma~\ref{lem:schauder}, with the truncated kernel \(c_{N,s}|z|^{-N-2s}\mathbf1_{\{|z|>\tau\}}\) and \(\tau\to0\) for fixed \(\delta>0\), gives
	\[
		(-\Delta)^s\phi_\delta=\rho_\delta*g\quad\text{pointwise in }\R^N.
	\]
	Thus \eqref{eq:Silvestre-holder} yields
	\[
		\|\phi_\delta\|_{C^{2s+\sigma}(\R^N)}
		\le C\bigl(\|\phi\|_{L^\infty(\R^N)}+\|g\|_{C^{0,\sigma}(\R^N)}\bigr).
	\]
	Arzel\`a--Ascoli and lower semicontinuity give a representative \(\phi\in C^{2s+\sigma}(\R^N)\). The same truncation argument gives \((-\Delta)^s\phi=g\) pointwise.

\section*{Acknowledgements}
	S.~Alarc\'on was supported by FONDECYT grants 1211766 and 1221365 and by Proyectos Internos USM PI\_LIR\_25\_11 and PI\_LIR\_26\_18. L.\,M.~Del Pezzo was partially supported by ANII (Agencia Nacional de Investigaci\'on e Innovaci\'on, Uruguay) grant FCE\_3\_2024\_1\_181302. A.~Quaas was partially supported by FONDECYT grant 1231585. L.\,M.~Del Pezzo and A.~Quaas were partially supported by grant PID2023-148028NB-I00, funded by MICIU/AEI/10.13039/501100011033 and by ERDF/EU.

\section*{Declarations}
	\noindent\textbf{Conflict of interest.} The authors declare that they have no conflict of interest.\medskip

	\noindent\textbf{Data availability.} Data sharing is not applicable to this article, as no datasets were generated or analyzed.

	\end{document}